%% file: recycle_imaging4.tex
\documentclass{siamltex}

\include{rom_macros}

\usepackage{subfigure}

\newcommand{\bfMM}{{\bf M}}

\newcommand{\V}[1]{{\boldsymbol{\mathbf{\MakeLowercase{#1}}}}} % vector
\newcommand{\M}[1]{\mathbf{\MakeUppercase{#1}}} % matrix
\newcommand{\Vv}{\V{v}}
\newcommand{\Vq}{\V{q}}
\newcommand{\Vr}{\V{r}}

\newcommand{\Vb}{\V{b}}

\newcommand{\MM}{\M{M}}
\newcommand{\MA}{\M{A}}

\newcommand{\MS}{\M{S}}
\newcommand{\MV}{\M{V}}
\newcommand{\MU}{\M{U}}
\newcommand{\MD}{\M{D}}

\newcommand{\MW}{\M{W}}
\newcommand{\ME}{\M{E}}

\newcommand{\MI}{\M{I}}
\newcommand{\MK}{\M{K}}
\newcommand{\MY}{\M{Y}}

\newcommand{\sI}{\mathcal{I}}

\title{Subspace Recycling for Sequences of Shifted Systems with Applications in Image Recovery \footnotemark[1]}

\author{ 
Misha E. Kilmer\footnotemark[3]\  and
Eric de Sturler\footnotemark[4]
}

\begin{document}
\maketitle
\renewcommand{\thefootnote}{\fnsymbol{footnote}}
\footnotetext[1]{This material is based upon work supported by the National Science
Foundation under Grants No. {NSF-DMS} 1720305,
1720291,
2208470, and {NSF HDR CCF}-1934553.}
\footnotetext[3]{Department of Mathematics, Tufts University, Medford, MA 02115.}
\footnotetext[4]{Department of Mathematics, Virginia Tech, Blacksburg, VA 24061.}
\renewcommand{\thefootnote}{\arabic{footnote}}

%%%%%%%%%%%%%%%ABSTRACT %%%%%%%%%%%%%%%%%
\begin{abstract}
For many applications involving a sequence of linear systems with slowly changing system matrices, subspace recycling,
which exploits relationships among systems and reuses search space information, 
can achieve huge gains in iterations across the total number of linear system solves in the sequence. 
However, for general (i.e., non-identity) shifted systems with the shift value varying over a wide range, the properties of the linear systems vary widely as well, which makes recycling less effective. If such a sequence of systems is embedded in a nonlinear iteration, the problem is compounded, and special approaches are needed to use recycling effectively.

In this paper, we develop new, more efficient, Krylov subspace recycling approaches for large-scale image reconstruction and restoration techniques that employ a nonlinear iteration to compute 
a suitable
regularization matrix. For each new regularization matrix, we need to solve regularized linear systems, $\MA + \c_\ell \ME_k$, for a sequence of regularization parameters, $\c_\ell$, to find the optimally regularized solution that, in turn, will be used to update the regularization matrix.   

In this paper, we analyze system and solution characteristics to choose appropriate techniques to solve each system rapidly.
Specifically, we use an inner-outer recycling approach
with a larger, principal recycle space for
each nonlinear step 
and smaller recycle spaces
for each shift.
We propose an efficient way to obtain good initial guesses from the principle recycle space and smaller shift-specific recycle spaces that lead to fast convergence.
Our method is substantially reduces the total number of matrix-vector products that would arise in a naive approach.   
Our approach is more generally applicable to sequences of shifted systems where the matrices in the sum are 
positive semi-definite.  

\end{abstract}

%%%%%%%%%%%%%%%%%%%% KEYWORDS %x%%%%%%%%%%%
%%%%%%%%%%%%%%%% AMS NUMBERS %%%%%%%%%%%%

%%%%%%%%%%%%%%% PAGE STYLE AND LABELING %%%%%%%%%%%%%
%\thispagestyle{plain}
\markboth{Kilmer, de Sturler, O'Connell}{authors}

\section{Introduction}
The need to solve sequences of large-scale, shifted, linear systems of the form 
\begin{equation}  \label{eq:seqsys} ( \bfA + \gamma_{\ell} \bfE_k ) \bfx^{(k,\ell)} = \bfb  \end{equation}
arises in many important applications, 
such as acoustics \cite{Biermann_2016,MeerBai_2010}, hyperspectral diffuse optical tomography \cite{ArvindhyDOT}, and electromagnetic problems
\cite{Haber_2015,ShaoPengLee_2011}. 
The expense of solving a large number of such systems is often the computational bottleneck of the larger (typically nonlinear) problem in which they are embedded. Several approaches have been developed to address the computational cost associated with solving sequences of shifted systems, though many focus on identity shifts 
\cite{FromGlas_1998,Sood_2016,Sood_2014}.
Some methods use Lanczos recurrences or Arnoldi iterations, see e.g., \cite{Darnelletal2008, KilmerdeSturler2006, Meerbergen2003, Saibabaetal2013, Simoncini2003}.  

The presence of a non-identity shift matrix complicates Krylov subspace recycling because the recycle space gets projected differently for each shift.  To use recycling efficiently over both shift $\ell$ and outer-iteration $k$ requires sophisticated new recycling tools. To this end, we propose several innovations.

In this paper, we develop new Krylov subspace recycling approaches for large-scale image reconstruction and restoration techniques that employ a nonlinear iteration to compute  
a suitable
regularization matrix such that edges are preserved \cite{Gazzola_2020,IRNekki,wohlberg2008lp,VogelOman,IRNtv}. Such computations are increasingly common in linear ill-posed problems where one has to learn the prior through a bootstrapping process.
These problems are expensive, as they involve three nested iterations. The outer iteration computes a sequence of improving regularization matrices or priors $\ME_k$; see below. For each regularization matrix, we solve regularized linear systems, $\MA + \c_\ell \ME_k$, for a sequence of regularization parameters, $\c_\ell$, to find the optimally regularized solution that, in turn, will be used to update the regularization matrix. Each large linear system is solved using an iterative method.    
As a result, the total number of inner-most iterations, counted as matrix-vector products, can be very large.

\comment{For this class of problems, Krylov subspace recycling can be highly effective in reducing the total number of matrix-vector products, both by providing fast convergence as well as good initial guesses. However, it is important to limit the overhead of potentially large numbers of repeated projections. To this end, we propose several innovations. }

\paragraph{\it Contributions}

First, we use the computations for the initial regularization matrix to seed a relatively large principal recycle space with both the solutions for a few selected shifts and the invariant subspace approximations constructed from the initial iterative solves. We show in Section \ref{sec:analysis} that seeding the space with the solutions for a few judiciously selected shifts results in good initial guesses for many shifts. Our proposed approach contrasts with other so-called ``seed methods" that use a single search space to accelerate the solution for multiple right hand sides \cite{ChanNg_1999, ChanWan_1997,KilMilRap_2001}.

Second, we use the principal recycle space to efficiently compute an optimal initial guess for each shift using small precomputed matrices that are independent of the shifts and that are computed once per outer iteration. The principal recycle space is updated for each outer iteration update of the regularization matrix. 
 
Third, if the computed initial guess is not sufficiently accurate, we use a recycling MINRES \cite{KilmerdeSturler2006,PaigeSaunders1975, Wangetal2007} with a relatively small shift-dependent recycling space for each system. 
After computing the block Rayleigh quotients of $\MA$ and $\ME$ for an orthogonal basis for the principal recycle space, approximate eigenpairs can be computed efficiently for each shift. In fact, the most expensive part of this computation can be reused from the work for computing the initial guesses.  These local recycle spaces are augmented with solutions that are nearby in shift and in outer iteration.

Fourth, we compute an effective update to the principal recycle space at the end of the outer iteration.

Since our systems, detailed in the next section, are representative of more general, square,  shifted systems (\ref{eq:seqsys}) 
where $\bfA$ and  $\bfE_k$ are symmetric and positive semidefinite,
this paper may be relevant to the larger class of applications mentioned previously.

The remainder of the  paper is organized as follows. 
In Section \ref{sec:back}, we introduce the necessary background, terminology, and notation for deblurring and edge preserving regularization, as well as subspace recycling. 
In Section \ref{sec:analysis}, we analyze the properties of the system matrices and solutions to set the stage for our algorithmic choices.
The main features of our new algorithms are given in Section \ref{sec:recycle}, where we introduce our inner-outer recycling approach for a fixed $\bfE_k$ over all the shifts.  
In Section \ref{sec:rspace}, we explain how to choose the various initial recycle spaces and how to expand these recycle subspaces based on information learned from the iterative-recycling solves. In this section, 
we also  give a pseudocode for our algorithm and 
some more detail and analysis. In Section \ref{sec:overhead}, we provide a detailed overview of the overhead of our proposed algorithmic changes and computational cost issues. Some alternative algorithmic choices and their pros and cons are discussed in Section \ref{sec:oblique}.
Section \ref{sec:num} is devoted to numerical results, and in Section \ref{sec:conclusions} we provide conclusions and future work.

\section{Background}  \label{sec:back}
We begin by reviewing features of the deblurring problems and necessary features of corresponding regularization tools that are used to recover edge information.  We then review one type of regularization parameter selection tool that we use in our experiments, called the L-curve approach (see, for example, \cite{hansenbook, 
LawsonHanson,DPOHansen} for details).  Our method is more or less independent of the parameter selection heuristic.   However, as the presentation of our approach is greatly improved by fixing the parameter selection algorithm and discussing the remaining features of our algorithm in the context of the selected method, we have chosen here to use one specific, well-known method, which does not require knowledge of the noise level.   
   
\subsection{Forward Operator and Regularization}
The class of systems we consider in this paper arises in linear inverse problems for the reconstruction or restoration of
images, 
where we need to learn the prior or regularization matrix 
to preserve edge features in 
the regularized solutions for those images. 
Suppose that $\bfx$ is a vectorized version of a 2D grayscale image, and let $\bfC$ denote the known linear forward operator (e.g., a blurring operator, Radon transform, etc.).   
The forward mapping from image space to data space is
\begin{equation}  \label{eq:forward}   \bfd = \bfC \bfx + \bfe, 
\end{equation}
where $\bfd$ denotes the measured data that has been corrupted by (unknown) noise $\bfe$.   Without loss of generality, we assume that
the noise, $\bfe$, is white\footnote{We assume that enough information is known to  whiten the data.}.

The typical approach is to solve 
\begin{equation} \label{eq:tikhonov} \min_{\bfx} \| \bfd - \bfC \bfx \|_2^2 + \lambda^2 \mathcal{R}[\bfx] , 
\end{equation}
where $\mathcal{R}[\bfx]$ denotes a regularization operator applied to $\bfx$.  
One popular choice to recover edges in the image is an approximation to total variation regularization (TV), in which $\mathcal{R}[\bfx] = \sum_i \sqrt{ (\bfDel_x \bfx)_i^2 + (\bfDel_y \bfx)_i^2 +\beta} $, where $\bfDel_x, \bfDel_y$ represent discrete derivatives in the x and y directions, $\beta$ is a small positive constant, and $i$ runs over all pixels. Alternative penalties include $l_p$-norm, $1 \le p <2$, penalties on the gradient image. 
Methods for solving (\ref{eq:tikhonov}) with such constraints are iterative in nature, requiring a sequence of linearly constrained systems to be solved (see \cite{IRNekki,wohlberg2008lp,Vogel,IRNtv}, for example).  
  
In this paper, we work with an alternative to TV introduced in \cite{Gazzola_2020}, which 
we discuss in some detail 
below. 
%in section \ref{sec:back}.  
% However, we note that this is
This method is representative of a larger class of regularizers that enforce edge preservation and the optimal $\lambda$ is not known a priori.  

Specifically, our shifted systems come in the form of regularized normal equations arising in the aforemented regularization:
\begin{equation} \label{eq:blkeqn} 
  \left( \bfC^T \bfC + \lambda_{\ell}^2 \bfL^T \bfD_k^2 \bfL \right) \Xkl = \bfC^T \bfd, 
\end{equation}
must be solved for $k=0,\ldots,K$ and $\ell = 1, \ldots,M$.

Consider again the forward model in (\ref{eq:forward}). The singular values of $\bfC$ decay over many orders of magnitude toward zero.   The large singular values are associated with the ``smooth'' modes, while the small singular values correspond to the high frequency modes.   Thus, the eigenvectors of 
$\bfA = \bfC^T\bfC$ associated with the large eigenvalues corresponds to smooth modes, and the eigenvectors associated with the small eigenvalues corresponds to high frequency.   Roughly the opposite is true with respect to $\bfE_k = \bfL^T \bfD_k^2 \bfL$, since it is related to a diffusion operator (see \cite{Gazzola_2020} and discussion and references therein).   Natural images, even those with much edge detail, contain large components in the directions of the smooth modes.    The trick is in recovering information associated with the edges which is buried in the subspace of $\bfA$ corresponding to the high frequency modes.

The method in \cite{Gazzola_2020} attempts to recover edge detail in images over a sequence of steps.  In the first step, the systems
\[  (\bfA + \gamma_{\ell} \bfE_0) \bfx^{(0,\ell)} = \bfb, \quad \mbox{ with } \quad \bfA:= \bfC^T \bfC, \qquad \bfE_0 = \bfL^T \bfL  \qquad \bfb = \bfC^T \bfd,\]
need to be solved (to reasonable approximation) for $\ell = 1,\ldots,M$.   The `optimal' value of $\gamma_{\ell}:= \lambda_{\ell}^2$ is selected to be the point on the parametric curve with respect to  $\lambda$ containing the discrete points $(\log(\| \bfC \bfx^{(0,\ell)} - \bfd \|),\log(\| \bfL \bfx^{(0,\ell)} \|)$ that corresponds to the corner of the curve (point of maximum curvature). 
This method of choosing the regularization parameter is called the L-curve method 
(see \cite{hansenbook} for details and references). We call the solution corresponding to the selected parameter $\bfx^{(0,*)}$.   

If the noise is not too large and the blur is not too great, then $\bfx^{(0,*)}$ will at least contain some transitional pixels that correspond to edges.   In other words, if we look at the ``gradient image'', $\V{q}^{(0)} := \bfL \bfx^{(0,*)}$, some edge information should be visible.  We do not want to penalize the edge information obtained in this first round in subsequent improvements of the image.   

The method works to improve the previous estimate as follows.  Define $\bfw = |\V{q}^{(0)}|/\| \V{q}^{(0} \|_{\infty}$, and then overwrite it by
$ \bfw \leftarrow 1- \bfw.^p$.  This ensures that real edges that have been picked up in the first stage have 0 or small weights associated with them.  Then, we set $\bfD_1 = \mbox{diag}(\bfw)$, and solve 
\[ 
(\bfA + \gamma_{\ell}  \bfE_1) \bfx^{(1,\ell)} = \bfb, \qquad \bfE_1 = \bfL^T \bfD_1^2 \bfL ,
\]
for all the $\gamma_\ell$ again.   We get a new L-curve:         $(\log(\| \bfC \bfx^{(1,\ell)} - \bfd \|),\log(\| \bfD_1 \bfL \bfx^{(1,\ell)} \|)$  since the regularization operator $\bfL$ has been replaced by $\bfD_1 \bfL$.  
The solution corresponding to the corner is called $\bfx^{(1,*)}$, and we expect it to have more edge information than the prior image.   

We repeat the process, but this time using the normalized gradient image: $\bfw = |\MD_1\V{q}^{(1)}|/\| \MD_1 \V{q}^{(1)}\|_{\infty}$, create a new weight matrix from $\bfw \leftarrow 1 - \bfw.^p$, so that $\bfD_{2} = \mbox{diag}(\bfw) \bfD_{1}$, and now repeat the 
solution process.   This continues until edges have been recovered as much as possible, as measured by tracking image gradient quality.  
     
The `optimal' regularization parameters are expected to be increasing in magnitude as a function of the outer iteration index $k$. See \cite{Gazzola_2020} for details.  
                       
\subsection{The L-curve, Shift Range, and Solver Accuracy}
   
For ease of discussion, we assume that $\| \bfC \|_2 \approx \| \bfL \|_2 \approx 1$.   In our application, we know that the entries in $\bfD_k$ are bounded by 1, and hence $\| \bfD_{k} \bfL \|_2 $ does not grow.   We do not know the optimal value of $\lambda$ for any fixed $k$ in advance.  In the L-curve method, the point of maximum curvature
of the curve
$\left( \log( \| \bfC \bfx^{(k,\ell)} - \bfd \|),\log( \| \bfD_k \bfL \bfx^{(k,\ell)} \|) \right)$
is chosen as the regularization parameter. 
For that point to be well-defined, however, the discrete values of $\lambda$ must span several orders of magnitude.  Moreover, the range is dictated by the relative sizes of $\| \bfC \|$ and $\| \bfL \|$, as well as the singular values of $ \bfC$ and the noise level.  
{\em Importantly, each of the systems needs to be solved sufficiently accurately for the regularizing effect of incorporating $\lambda$ to dominate 
the known regularizing effects of the iterative solver.}
If the systems are not solved accurately enough, particularly for the smallest $\lambda$ values, we see only the regularizing effects of the iterative solver and 
the points on the L-curve appear so clustered that
the corner is poorly defined.
Likewise, if 
we do not choose large enough values, we may not see the tail of the L-curve. This is why, in the numerical results section, we solve each system to the same relative residual norm tolerance. 

Given the assumed normalization of the problem, we will use M, log-evenly spaced regularization parameters between $10^1$ and $10^{-4}$ in our experiments, meaning the shifts $\gamma_\ell = \lambda^2_{\ell}$ differ by 10 orders of magnitude.  Given our initial normalizations
on $\bfC$ and $\bfL$, and the noise levels used, this range is sufficient to ensure that if all systems are solved to the same relative residual value, we can discern the corner of an L-curve.   Note that because our penalty
term is $\bfD_{k} \bfL$, the meaning of the
vertical axis of the L-curve will change at each outer iteration $k$.

\subsection{General Recycling}
We first discuss recycling for a single system; see \cite{KilmerdeSturler2006}.
More on recycling can be found in the literature; see, for example, \cite{Ahujaetal2012, KilmerdeSturler2006, Parksetal2006, Wangetal2007,SoodStuKil2020}.

Consider the linear system $\bfA \bfx=\bfb$, with symmetric
$\bfA \in \mathbb{R}^{N \times N}$ and $\bfb \in \mathbb{R}^N$.    
Given $\widetilde{\bfUU} \in \mathbb{R}^{N \times n_c}$, we 
compute an initial
approximate solution over $\mrg(\widetilde{\bfUU})$ using the
skinny
QR factorization $\bfK \bfR = \bfA \widetilde{\bfUU}$  
and 
$\bfUU = \widetilde{\bfUU} \bfR^{-1}$.
Minimizing the residual norm, $\| \bfb  - \bfA \bfx \|_2$, under the constraint $\bfx = \bfUU \bfz$, for some $\bfz$,
gives the Petrov-Galerkin condition 
$\bfK^T(\bfb - \bfA \bfUU\bfz) = 0$, 
and the 
approximate solution
\begin{equation}
\bfx_0 = \bfUU \bfK^T \bfb,
\label{eq:initX}
\end{equation}
with initial residual $\bfr_0 = (\bfI-\bfK \bfK^T) \bfb$.  If the relative residual norm is not sufficiently small, we expand the space in which we look for a solution
using the Lanczos recurrence with $(\bfI-\bfK \bfK^T)\bfA$ and  $\bfv_1 = (\bfI-\bfK \bfK^T)\bfb/\|(\bfI-\bfK \bfK^T)\bfb\|_2$.  This gives
\eqs
  (\bfI-\bfK \bfK^T) \bfA \bfV_m & = & \bfV_{m+1}\underline{\M{T}}_m 
  \;\; \Leftrightarrow \;\;  
  \bfA \bfV_m = \bfK \bfK^T \bfA \bfV_m + \bfV_{m+1}\underline{\M{T}}_m ,
\label{eq:eq01}
\eqe
where 
$\underline{ \M{T} }_m$ is an 
$(m+1) \times m$ tridiagonal matrix.\footnote{A short recurrence is possible for $(\bfI-\bfK \bfK^T) \bfA$, as the matrix is self-adjoint (symmetric) over the Krylov space generated from starting vector $\bfv_1$ \cite{KilmerdeSturler2006,SoodStuKil2020}.}

Next, we consider the expanded space $\mrg ([\bfUU \; \bfV_m ])$ and compute the approximate solution
that minimizes $\|\bfb - \bfA( \bfV_m \bfy + \bfUU \bfz)\|_2$ as follows:
\eqs
  \min_{\bfy,\bfz} \left\| \bfb -\bfA[\bfUU \; \bfV_m ]
  \left[\begin{array}{c}
    \bfz \\ \bfy
  \end{array}\right]
  \right\|_2 &  &
\nonumber \\
  = \min_{\bfy,\bfz} \left\| \bfb - [ \bfK \; \bfV_{m+1}]
  \left[\begin{array}{cc}
    \bfI & \bfK^T\bfA \bfV_m  \\
    0 & \underline{\M{T}}_m\\
  \end{array}\right]
  \left[\begin{array}{c}
    \bfz \\ \bfy
  \end{array}\right] \right\|_2 
  &  &
\nonumber \\
\label{eq:SmallMin}
  = \min_{\bfy,\bfz} 
    \left\| \left[\begin{array}{c}
      \Co^T \bfb \\ \xi \bfe_1
    \end{array}\right] -
    \left[\begin{array}{cc}
      \bfI & \bfK^T\bfA \bfV_m     \\
      0 & \underline{\M{T}}_m\\
    \end{array}\right]
    \left[\begin{array}{c}
      \bfz \\ \bfy
    \end{array}\right] \right\|_2 ,
\label{eq:eqrecyc}
\eqe
where $\bfe_1$ denotes the first Cartesian basis vector in $\mathbb{R}^{m+1}$
and $\xi = \|(\bfI - \bfK \bfK^T)\bfb \|_2$.
We first compute the solution to the projected problem
\[ 
  \min_{\bfy} \| \underline{\M{T}}_m \bfy - \xi \bfe_1 \|_2, 
\]
then compute 
$\bfz = \bfK^T\bfb -  \bfK^T\bfA\bfV_m \bfy$ to satisfy
(\ref{eq:SmallMin}), and finally
set $\bfx = \bfV_m \bfy + \bfUU \bfz$. 
The solution component $\bfV_m \bfy$ can be computed via short term recurrences \cite{Melloetal2010,Wangetal2007} as for standard MINRES, so 
we do not have to explicitly store $\bfV_m$.

Finally, if $n_c$ is large, the additional orthogonalizations in each iteration form a source of overhead that cannot be ignored,
hence it is important to keep the number of columns in $\bfK$, $n_c$, relatively small.

\subsection{Shifted Systems}
It is easy to see that adding a shift poses a problem with recycling. We have 
$\bfA \bfUU = \bfK$,  $\bfK^T \bfK = \bfI$, but $\mrg(\bfA + \gamma_\ell \bfE) \bfUU$ will not coincide with
$\mrg(\bfK)$ in general, and will 
change with $\gamma_\ell$.
Various solutions have been proposed in 
\cite{KilmerdeSturler2006,OConKilStuGug17, Sood_2016}; see also
\cite{BauGijz_2015,FromGlas_1998,SoodStuKil2020,Sood_2014}.

In \cite{OConKilStuGug17}, the authors propose an inner-outer Krylov recycling scheme, primarily to address a sequence of
systems of the form
\[ \bfA^{(k)} \bfx^{(k)}_j = \bfb_j, \qquad j=1,\ldots,n_{s} .\]
The idea was to maintain a {\it single, larger, principle recycle space} in which to look for initial solutions to all the right-hand sides for a fixed $k$ - this is the outer part of the scheme.  Then, {\it local recycle spaces} were constructed from small subsets of the columns of this principle space.  
If, for system $j$, the relative residual norm was not small enough, they solved a new system for an updated guess. The updated guess was constrained to be in 
the sum of the current local recycle space and the Krylov subspace with right-hand side generated from the current residual.   The principle and local recycle spaces {were only expanded by the update to the solution estimate that was not already contained in the principle and local recycle spaces.}  
In this way, smaller recycle spaces were used for actual iterative solves, keeping the overhead involved with the recycling to a minimum.  

We will adopt a similar strategy in the present paper: 

\medskip
\begin{center}\fbox{ \parbox{4in}{We establish a larger, principle recycle space, for determining\\ initial guesses,
and smaller, tailored, local subspaces for recycling\\ for correction equations on individual systems.}}   \end{center}

\medskip
However, the method for determining the local spaces and the updates is decidedly different from previous approaches.  Before detailing our strategy in Sections \ref{sec:recycle} and \ref{sec:rspace}, we provide an analysis in the next section that
will guide us in the process of defining and updating recycle spaces.

\section{Analysis}  \label{sec:analysis}
In this section, we analyze the properties of the system matrices in (\ref{eq:seqsys}) and the 
solutions to (\ref{eq:seqsys}) to guide our 
strategies for efficient solution.

We assume that $\bfA$ is SPD, and 
$\bfE$ is SPSD.  
The analysis that follows requires only those two conditions.
In our particular application, since $\bfA := \bfC^T \bfC$, we can ensure that $\bfA$ is SPD by assuming that $\bfC$ has full column rank.  Likewise, $\bfE$ is SPSD by virtue of its definition in
our application.  However, the theory below is applicable
for other applications for which $\bfA$ is SPD and factors as
$\bfC^T \bfC$, and for which $\bfE$ is
SPSD.

\subsection{Generalized Eigenvector Analysis}
First we analyze the solutions  to (\ref{eq:seqsys}) and how they vary as a function of the shift $\c_j$.

Given our assumptions, there exists $\bfV$ that solves the generalized eigenvalue problem
  \[ \bfE_k \bfV = \bfA \bfV \bfMM, \qquad \bfMM = \mbox{diag}(\mu_1,\ldots,\mu_n), \]
so that $\bfV^T \bfA \bfV = \bfI$, and hence the vectors $\bfC \bfv_j$ are an orthonormal set.  
This also implies that $\MV^T \ME_k \MV = \MM$.
We further assume that $\mu_1 \ge \mu_2 \ge \cdots$.   
There must be at least 1 zero eigenvalue ($\mu_n = 0$) in our application, since
$\bfL$ has a non-trivial null space. 
Consider the solution in the $\bfV$ basis as a function of $\gamma$: $\bfx_\gamma = \bfV \bff_\gamma$ for some $\bff_\gamma$. Then 
\eqsn  
  \bfA \bfV \bff_\gamma + \gamma \bfE_k \bfV \bff_\gamma = \bfC^T \bfd \quad 
  & \eqv & \quad \bfV^T \bfA \bfV \bff_\gamma + \gamma \bfV^T \bfE_k \bfV \bff_\gamma  = \bfV^T \bfC^T \bfd \eqv \\ 
  (\bfI + \gamma \bfMM) \bff_\gamma = \tilde{\bfd}
  \quad & \eqv & \quad 
  \bff_{\c} = 
  \left( 
    \frac{\wtl{d}_1}{1+\c \m_1} \;\; 
    \frac{\wtl{d}_2}{1+\c \m_2} \;\; 
    \cdots \;\;
    \frac{\wtl{d}_n}{1+\c \m_n}
  \right)^T ,
\eqen
where $\tilde{\bfd}$ is the vector of expansion coefficients of $\bfd$ in the orthonormal basis given by the columns of the matrix $\bfC \bfV$.
For the solution $\bfx_\c$ this gives
\eqs\label{eq:x_gamma}
  \bfx_\gamma & = &
  \sum_{j=1}^N \bfv_j \frac{\wtl{d}_j}{1+\c \m_j} ,
\eqe
which provides two important observations. First, if $\c \m_j \ll 1$, then the coefficient along $\bfv_j$ is approximately $\wtl{d}_j$, and this holds for all smaller $\c$ as well. So, these solution components are insensitive to $\c$, and $\bfx_\c$ changes little in these directions until $\c \m_j \approx 1$. For example, if $\c\m_s \leq 1/100$, for all smaller $\c$ the change along $\Vv_s$ is less than $(1/100) \wtl{d}_s$.
Second, if $\c \m_j \gg 1$, then the coefficient along $\bfv_j$ is approximately $0$, and this holds for all larger $\c$ as well.
Again, these solution components are insensitive to (larger) $\c$, and $\bfx_\c$ changes little in these directions until $\c \m_j \approx 1$. Analogous, to the example above, if $\c\m_{\ell} \geq 100$, the change along $\Vv_{\ell}$ for an increase in $\c$
is bounded by $(1/100) \wtl{d}_{\ell}$.
As a result, for a change in $\c$ of, say, a factor $10$ (increase) only a modest number of components have an appreciable change. 
For any pair of shifts $\gamma_a, \gamma_b$ (WLOG, assume $\gamma_b > \gamma_a$), we have
\begin{equation} \label{eq:compare} 
(\bfx_{\gamma_a} - \bfx_{\gamma_b}) = \bfV(\bff_{\gamma_1} - \bff_{\gamma_2}) = \sum_{i=1}^{N-1}  \left( \frac{(\gamma_b - \gamma_a)\mu_i}{(1+\gamma_b \mu_i)(1+\gamma_a \mu_i)} \right) \tilde{d}_i \bfv_i .
\end{equation}

The analysis suggests that seeding a recycle space with a few solutions from the range of interest of the $\c$-values should be a good idea.  It also suggests that we may wish to work differently
across the spectrum of $\gamma$ values, keeping this relative gap, $(\c_b - \c_a)\m_i / (1+\c_b\m_i)$, small.  We will exploit both of these ideas in our
algorithm development.

\section{Two-level Recycling for Shifted Systems} \label{sec:recycle}

Consider again the sequences 
\begin{equation} \label{eq:shift} 
\left( \bfA + \gamma_\ell \bfE_k \right) \bfx^{(k,\ell)}= \bfb, \qquad \ell = 1,\ldots,M; \qquad k=0,\ldots,K ,
\end{equation} 
where
$\bfA = \bfC^T \bfC$, $\bfE_k = \bfL^T \bfD_k^2 \bfL$, and $\bfb = \bfC^T \bfd$ and $\gamma_\ell = \lambda_{\ell}^2$.   
 
As alluded to previously, our strategy is to build and maintain one larger, principle recycle space from which the initial solutions across all shifts will be determined. 
Then, we identify and/or update smaller, shift-specific (a.k.a local) recycle spaces from which to compute the corrections to the solutions for each system.   As mentioned earlier, we do not know a priori the value of the appropriate $\lambda$ for any given $k$.  
To get a good estimate for $\lambda$ via the L-curve, we need to solve the systems (\ref{eq:shift}) accurately enough that the L-curve, as a function of $\lambda$, has a well defined corner; see section \ref{ssec:resi}.  
  
In this section, we describe the two-level recycling {\em for a fixed $k$} and with the assumption that the principle and local recycling spaces have been determined.   
In the next section, we describe how to build the spaces initially, how to update the search spaces, and how to adapt our basic strategy for applications, such as ours, where the shifts vary greatly.

\subsection{Initial Guesses} \label{ssec:initorth}

Assume that $\widetilde{\bfUU} \in \mathbb{R}^{N \times n_c}$ contains a suitable principle recycle space.  
 Here, although $n_c \ll N$, $n_c$ is assumed to be large enough that using the principle recycle space is too expensive for recycling in the iterations.  However, $n_c$ is small enough that we can use the whole space for computing good initial guesses at a relatively low cost.  We 
 determine the initial guesses based on orthogonal projection. 
 Under the right circumstances, e.g., with a smaller shift range,
 oblique projections may be cheaper while providing almost equivalent results. This is not the case in our application.  Therefore, we defer the description -- included for the sake of completeness -- to Section \ref{sec:oblique}.

 Recall that $0 < \gamma_1 < \gamma_2 < \cdots < \gamma_M$.   
Let $\bfUU \in \mathbb{R}^{N \times n_c}$ be generated from the recycle space $\widetilde{\bfUU}$ such that
\begin{equation} \label{eq:establishU} 
  \left( \bfA + \gamma_* \bfE_k \right) \bfUU = \bfK, \quad\mbox{and}\quad \bfK^T \bfK = \bfI , 
\end{equation}
where for the purpose of this discussion $\c_*$ is
an arbitrary choice.
To compute an initial solution from the principle recycle space, $\mrg (\bfUU)$, 
for each $\c_\ell$, 
we set $\Xkl_0 = \bfUU \bfq$,
and solve for $\bfq$.
% (for some $\bfq$). 
The initial residual is given by
\begin{eqnarray} \label{eq:res_initGuessPrincRecy}
\nonumber
  \bfrkl_0  & = & \bfb - \left( \bfA + \gamma_\ell \bfE_k \right) \bfUU \bfq
   = \bfb - \left( \bfA + \gamma_* \bfE_k + (\gamma_\ell - \gamma_*) \bfE_k \right) \bfUU \bfq 
\\ 
   & = & 
   \bfb - \left( \bfK  + \d_\ell \bfE_k \bfUU \right)\bfq,
\end{eqnarray}
where $\d_\ell = \gamma_\ell - \gamma_*$.
To minimize $\| \Vr_0^{(k,\ell)} \|_2$, 
we pick $\Vq$ such that $\Vr_0^{(k,\ell)} \perp 
\bfK  + \d_\ell \bfE_k \bfUU$, which gives (using the 
symmetry of $\ME_k$)
\eqs 
\nonumber
  (\bfK  + \d_\ell \bfE_k \bfUU)^T \Vr_0^{(k,\ell)}
  & = & 
  (\bfK  + \d_\ell \bfE_k \bfUU)^T
 (\bfb - \left( \bfK  + \d_\ell \bfE_k \bfUU \right)\bfq) = \V0 \quad \eqv \\
 \MK^T\Vb + \d_\ell \MU^T\ME_k\Vb
 & = & 
 (\MI + \d_\ell(\MK^T\ME_k \MU) + 
 \d_\ell(\MK^T\ME_k \MU)^T + 
 \d_\ell^2 \MU^T\ME_k^2\MU) \Vq, \label{eq:orthupdate}
\eqe
where only $\d_\ell$ changes, and the $n_c \times n_c$
matrices $\MK^T\ME_k \MU$ and 
$(\ME_k\MU)^T(\ME_k\MU)$ and the vector
$\MU^T\ME_k\Vb$ have to be computed only once
(per outer iteration $k$).
Computing a Cholesky factorization for each
$\d_\ell$ and solving for $\Vq$
has only a $O(n_c^3)$ cost and does not
involve any $O(N)$ work.  

The discussion of the overhead in this step as well as all other parts of the algorithm will be deferred until Section \ref{sec:overhead}, after all algorithmic components have been presented, and all overheads can be assessed together.

\comment{\medskip
\noindent{\bf Cost Summary:} the $O(N)$ or greater effort involved with this step is:
\begin{enumerate}
    \item $n_c$ matvecs with $\MA, \ME$ in (\ref{eq:establishU}).
    \item Orthogonalization of $n_c$ vectors to get $\MK$ in (\ref{eq:establishU}). Implicitly apply $\M{R}^{-1}$ to the left of $\tilde{\MU}$ and to left of the precomputed $\ME_k \tilde{\MU}$ to obtain $\ME_k \MU$. ($2N$ lower triangular $n_c \times n_c$ solves.) 
    \item Form $\MK^T (\ME_k \MU)$, $(\ME_k \MU)^T(\ME_k \MU)$, and $(\MU_k^T\ME_k)\V{b}$, 
    using previously computed $\ME_k \MU$. ($O(n_c^2 N)$ flops).
\end{enumerate}
}

\subsection{Shift Specific Recycling}
With all initial guesses and residuals computed, we check the relative residual
norms against our convergence criterion. 
If an initial residual is not suitably small, we solve only for the required  
%the incremental change from this 
update to the initial guess, 
$\bfg^{(k,\ell)} = \Xkl - \Xkl_0$.
%, 
%the goal is to recover . 
The correction equations are  
\begin{eqnarray} \label{eq:corr_g}
\left( \bfA + \gamma_\ell \bfE_k \right) \bfg^{(k,\ell)}  =  
%\underbrace{
\bfb - \left( \bfA + \gamma_\ell \bfE_k \right) \Xkl_0 = \bfr^{(k,\ell)} ,
\quad 
\mbox{for }
\ell = 1, \ldots, M.
%_{\bfr^{(k,\ell)}} ,
\end{eqnarray}

Suppose for the moment we have a skinny local recycle matrix $\widetilde{\bfUU}_\ell$ such that
   \begin{equation} \label{eq:localU} (\bfA + \gamma_\ell \bfE_k) \widetilde{\bfUU}_\ell = \bfK_\ell \bfR_{\ell}, \qquad \bfK_\ell^T \bfK_\ell = \bfI, \qquad \bfUU_{\ell} := \widetilde{\bfUU}_{\ell} \bfR_{\ell}^{-1} ,  \end{equation}  
The subscript indicates that this is a {\bf shift specific recycle space}.  We address how to choose this space 
in Section \ref{ssec:defineSS}.
Given this space, 
to get $\bfg^{(k,\ell)}$, we solve
\[\min_{\bfg \in \mathcal{S}_\ell} \| \bfrkl - \left( \bfA + \gamma_\ell \bfE_k \right) \bfg \|_2,
\]
over an appropriate {\bf local} subspace $\mathcal{S}_\ell$.  
We choose $\mathcal{S}_\ell = \Ra{[\bfUU_\ell,\bfV_{m_\ell}]}$, where $\Ra{\bfUU_\ell}$ is 
the local recycle space, and the columns of $\bfV_{m_\ell}$ form the basis for the $m_{\ell}$-dimensional Krylov subspace  
generated by the matrix $(\bfI-\bfK_\ell \bfK_\ell^T)\left( \bfA + \gamma_\ell \bfE_k \right)$ and 
\[
\bfv_1 = \left( \left( \bfI - \bfK_\ell \bfK_\ell^T \right)\bfrkl \right) /\|\left( \bfI - \bfK_\ell \bfK_\ell^T \right)\bfrkl\|_2.
\]
The subscript $\ell$ on the iteration count $m_{\ell}$ indicates that the number of iterations
to reach
convergence varies 
from one system to the next.

\paragraph{Remark} We 
emphasize that the present approach  
differs from the 
inner-outer recycling 
in \cite{OConKilStuGug17}. 
Here, $\bfK_{\ell}$
varies with the shift and the shift specific recycle space, whereas in 
\cite{OConKilStuGug17} 
$\bfK_\ell$ 
is always computed using $\bfA$
and varies only with the chosen right-hand side
dependent recycle space.  
Moreover, 
the initial solution estimates,
discussed above, are computed 
differently than in \cite{OConKilStuGug17}. 
 
Letting $m := m_{\ell}$ for ease of notation, using the 
Lanczos recurrence we obtain
\eqs
  (\bfI-\bfK_\ell \bfK_\ell^T) \left( \bfA + \gamma_\ell \bfE_k \right) \bfV_m & = & \bfV_{m+1} \underline{\M{T}}_m \Leftrightarrow \nonumber \\
  \left( \bfA + \gamma_\ell \bfE_k \right) \bfV_m & = & \bfK_\ell \bfK_\ell^T \left( \bfA + \gamma_\ell \bfE_k \right) \bfV_m + \bfV_{m+1} \underline{\M{T}}_m.  
\nonumber
\eqe
We then find $\bfy, \bfz$ by solving
\begin{equation}
\min_{\bfy,\bfz} \left\| \bea{c} \bfK_\ell^T \bfrkl \\ \xi \bfe_1 \eea -
      \bea{cc} \bfI & \bfK_\ell^T \left( \bfA + \gamma_\ell \bfE_k \right) \bfV_m \\
      0 & \underline{\M{T}}_m
      \eea \bea{c} \bfz \\ \bfy \eea \right\|_2 , 
      \label{eq:minShift}
\end{equation}
where $\xi = \|\left( \bfI - \bfK_\ell \bfK_\ell^T \right)\bfrkl\|_2$; see (\ref{eq:eqrecyc}) and the subsequent
discussion.
%Therefore, 
This gives the correction  $\bfg^{(k,\ell)} = \bfV_m\bfy +
\bfUU_\ell \bfz = 
\bfy_m + \bfUU_\ell(\bfK_\ell^T \bfrkl -\bfK_\ell^T \left( \bfA + \gamma_\ell \bfE_k \right) \bfV_m \bfy)$, and
\begin{equation}
\Xkl = \bfy_m + \bfUU_\ell \bfK_\ell^T [ \bfrkl -  \left( \bfA + \gamma_\ell \bfE_k \right) \bfy_m] + \Xkl_0. 
\label{eq:xFinal} 
\end{equation}

\section{Determining Recycle Spaces} \label{sec:rspace}

We now discuss how to initialize and update the principle recycle space  and 
the shift specific recycle spaces.

\subsection{Initializing the Principle Recycle Space}

Our principle recycle space needs to be selected such that we can obtain good initial solutions across all the $\gamma_\ell$.   
The properties of the system matrix $(\bfA + \gamma_\ell \bfE_k)$ change with $\gamma_\ell$, as do the properties of the solution.  When $k=0$, following the 
analysis in Section \ref{sec:analysis}, the systems for the smallest shifts are ill-conditioned 
as the eigenvalues of $\MA$ dominate.   
For the 
smallest shifts, this also means the solutions 
(which will be noise contaminated) 
are close, with a lot of high frequency information.

As the shifts increase, 
$\gamma_\ell \ME_0$ becomes a more dominant component of the operator. Hence, the eigenspace corresponding to the smallest eigenvalues of 
$(\bfA + \gamma_\ell \bfE_0)$
can no longer be categorized as completely high frequency.  
We need to account for this qualitative change in the solution spaces.  

Thus, for $k=0$, we initialize the principle recycle space as follows. % $\widetilde{\bfUU}$, 
First, using MINRES, we solve two\footnote{One could use more, as the application and shift spacing dictates.} representative systems, 
say, for $\gamma_{i_1}$ and $\gamma_{i_2}$, to residual convergence tolerance.  
With these two solves:
\begin{enumerate}
  \item We compute solution estimates $\bfx^{(0,i_1)}$ and  $\bfx^{(0,i_2)}$. 
  \item We get Ritz estimates for an approximate invariant subspace for the $\gamma_{i_1}$ (smaller shift) system.  
  Let the columns of $\MV_{i_1}$ span this space.      
  \item We also get Ritz estimates for an approximate invariant subspace for the $\gamma_{i_2}$ (larger shift) system.  
  Let the columns of  $\bfV_{i_2}$ span this space.   
\end{enumerate}
So, we initially set 
\begin{equation} \label{eq:widetildeU} \widetilde{\bfUU} = [ \bfx^{(0,i_1)},\bfx^{(0,i_2)},\bfV_{i_1},\bfV_{i_2}].
\end{equation}
As this matrix could be numerically rank deficient, 
we 
refine $\widetilde{\bfUU}$ as
in Algorithm \ref{alg:orthU}.
\begin{algorithm}[ht]  \caption{\label{alg:orthU} Stabilized $\widetilde{\bfUU}$}
\SetAlgoLined
   Compute skinny SVD: $\widetilde{\bfUU} = \M{Q} \M{S} \M{W}^T$ \\
   Determine $n_c$ such that $\frac{\sigma_1}{\sigma_i} < 10^{10}$ for all
      $1 \le i \le n_c$. \\
   Set $\widetilde{\bfUU} \leftarrow \M{Q}_{:,1:n_c}$\\
\end{algorithm}

\noindent
Although this involves some overhead, 
having a $\widetilde{\bfUU}$ with orthogonal columns is useful in computing shift-specific recycle spaces, as we will see in Section \ref{ssec:defineSS}.
So, this cost is amortized across all the computations.

\subsection{Updating the Principle Recycle Space for $k \ge 1$} \label{ssec:updatemaster}
The larger, principle recycle space serves
two purposes. 

First, 
at the start of the main loop, it is used to generate initial solutions via 
(\ref{eq:orthupdate}).
For $k \ge 1$, we have solutions $\bfx^{(k-1,\ell)}$ for $\ell = 1,\ldots,M$, and  
our analysis in Section \ref{sec:analysis} suggests that the solutions are related, particularly for the systems with small shifts.  Therefore, we include all $\bfx^{(k-1,\ell)}$ in the principle space
for the next outer iteration.
 
Second, the principle recycle space is also used to compute the Ritz approximations for generating the local recycle spaces used in 
the iterative solves (see the next subsection).  
For recycling to be effective, 
it should remove eigenvalues from 
$(\MA + \gamma_{\ell} \ME_k)$
that slow down 
convergence and remove large
spectral components of 
$\V{r}^{(k,\ell)}$.

The approximate invariant subspace $\mrg(\bfV_{i_1})$ corresponds to the small eigenvalues
of $\MA + \gamma_{i_1} \ME_0$.  We expect $\bfV_{i_1}$ to provide 
a good approximation to invariant subspaces of $\MA + \gamma_\ell \ME_k$ for the smallest shifts for all $k>0$.  
Therefore it is advantageous to keep $\bfV_{i_1}$ in our recycle space even for $k > 0$. 
On the other hand, as the outer iteration $k$ proceeds, 
the spectrum of $\bfE_k$ changes
due to better localization of the edges in the image. Since the eigenstructure 
of 
$\MA + \gamma_\ell \ME_k$ for large shifts is dominated
by $\ME_k$, $\bfV_{i_2}$ is not likely 
to be an adequate choice to retain in the principle space.

For a better alternative, we add and subtract $\gamma_\ell \bfE_{k-1}$ from (\ref{eq:seqsys}) to get  
\begin{equation} \label{eq:delta} \left(\frac{1}{\gamma_\ell}\bfA + \bfE_k \right)\delta \bfx^{(k,\ell)} = (\bfE_{k} - \bfE_{k-1})\bfx^{(k-1,\ell)} , 
\end{equation}
where $\delta \bfx^{(k,\ell)} := 
\bfx^{(k-1,\ell)} - \bfx^{(k,\ell)}$.  
We solve this equation for one of the large $\gamma_\ell$ (large $\ell$) values.
With $\bfx^{(k-1,\ell)}$ in the principle recycle space, 
adding an approximation to 
$\delta \bfx^{(k,\ell)}$ 
will produce 
a good initial guess $\bfx^{(k,\ell)}_0$. 
Hence, we take a few MINRES steps solving(\ref{eq:delta}). 
This also produces additional Krylov vectors, denoted $\MV^{(new)}$, that 
can be used for improved Ritz approximations to eigendirections that we want to project out
and for initial solutions for nearby systems.  
Therefore, we replace $\MV_{i_2}$ with 
the estimated solution correction and
a few Krylov vectors from solving this correction equation for one choice of $\ell$, say $\ell=\ell_c$.  
We fix the max number of iterates at 100 in our numerical examples, so that $\MV^{(new)}$ has at most 101 columns. 

Thus, our updated principle recycle space for $k > 0$ is initialized as
\begin{eqnarray} \label{eq:tildeUupdate} \widetilde{\bfUU} = [\MV_{i_1},\MV^{(new)},\bfx^{(k-1,1)},\ldots,\bfx^{(k-1,M)}], 
\quad
 \widetilde{\bfUU} \leftarrow \mbox{\tt stabilize}(\widetilde{\bfUU}),
\end{eqnarray}
where $\tt stabilize()$ refers to Algorithm \ref{alg:orthU}.

\subsection{Defining Shift-specific Recycle Spaces}
\label{ssec:defineSS}
We need shift specific recycle spaces such that convergence is fast for the correction equations (\ref{eq:corr_g}).  The operators are SPD, so typically, if we can remove directions corresponding to small eigenvalues, we can speed up convergence. 
We leverage the
principle recycle space, spanned by the orthonormal columns of $\widetilde{\bfUU}$, by 
efficiently
projecting one or more shifted systems
onto this space and 
initialize
the local recycle space using
selected Ritz vectors.

We assume that
\begin{equation} \label{eq:hq} \hat{\bfA} := \widetilde{\bfUU}^T \bfA \widetilde{\bfUU} \qquad \mbox{and} \qquad  \hat{\bfE}_k := \widetilde{\bfUU}^T \bfE_{k} \widetilde{\bfUU}, \end{equation}
have been precomputed. 
We select a shift $\gamma$, compute the projected shifted system, and identify Ritz vectors corresponding to the smallest Ritz values.  These are returned in a matrix $\MW$.
The process is described in Algorithm \ref{alg:interior}. 
\begin{algorithm}[ht]  \caption{\label{alg:interior} Generate Ritz Recycle Information for fixed $\gamma,\hat{\bfA},\hat{\bfE}$}
\SetAlgoLined
    Set $\bfH= \hat{\bfA} + \gamma \hat{\bfE}_k$, with $\hat{\bfA},\hat{\bfE}_k$ from  (\ref{eq:hq}) \\
   Find the eigenvectors, $\widetilde{\bfW}$, of this $n_c \times n_c$ matrix. \\
   Let $\mathcal{J}$ denote the indices corresponding to the smallest eigenvalues.\\
   Define $ \bfW = \widetilde{\bfUU}\widetilde{\bfW}(:,\mathcal{J})$ \\
\end{algorithm}

\begin{figure}
\centering
\includegraphics[scale=.4]{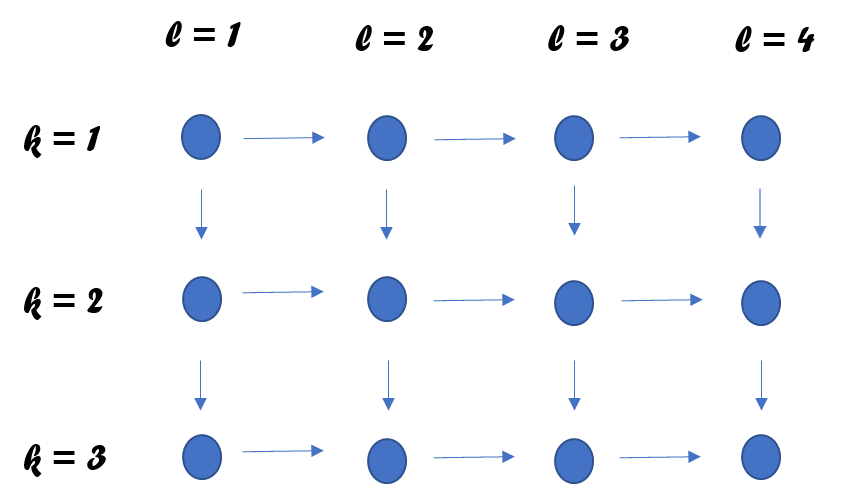}
\caption{\label{fig:local} Depiction of which $\bfg^{(k,\ell)}$ will be propagated forward and incorporated into the current local recycle space via Algorithm \ref{alg:updatelocal}. }
\end{figure}

For fixed $k$ and $\ell > 1$, we have just computed $\bfg^{(k,\ell-1)}$, and  for $k \ge 1$ and fixed $\ell$, we have the previous correction vectors $\bfg^{(k-1,\ell)}$. Figure \ref{fig:local} illustrates 
the information that can be propagated forward to update the local recycle spaces.\footnote{Since it is possible that $\bfg^{(k-1,\ell)}$ and $\bfg^{(k,\ell-1)}$ contain nearly redundant information, in practice we check that the angle between them is not close to 0, at the cost of one inner product. If it is, we add only one of them.}
%\begin{center}
%\end{center}
%  
The procedure to update the local spaces $\bfUU_{\ell}$ is
given in Algorithm \ref{alg:updatelocal}.  In the last step of Algorithm \ref{alg:updatelocal}, if $k=0$ and/or $\ell = 1$ so that either of the two right-most vectors in the update is undefined, it is understood they are omitted.   

\begin{algorithm}[ht] \caption{ \label{alg:updatelocal} 
Update local recycle space $\bfUU_{\ell}$}
 If $\MW$ not passed in, call Algorithm \ref{alg:interior} with suitable $\gamma$ and return $\MW$ \\
 Update $\widetilde{\bfUU}_{\ell}$ as $[\MW,\bfg^{(k-1,\ell)},\bfg^{(k,\ell-1)}]$ \\
 \end{algorithm}   

Now $\MW$ has $\le | \mathcal{J} |$ columns, and this is much smaller ($\approx 10$) than the dimension, $n_c$, of the principle recycle subspace.  Thus $\widetilde{\M{U}}_\ell$ has no more than about a dozen columns,
resulting in very modest overhead (see Section \ref{sec:overhead}).

\comment{
\medskip
\noindent{\bf Cost Summary}.  The $O(N)$ or greater costs are:
\begin{enumerate}
 \item The overhead in determining the local recycle spaces consists of\footnote{If $\widetilde{\bfUU}$ does not already have orthonormal columns, we need to include the cost to compute the orthonormal basis $\bfQ$.} computing each of the projections $\hat{\bfA}, \hat{\bfE}_k$.  This {\it only needs to be done only once per outer iteration $k$, independent of how many calls to Algorithm \ref{alg:interior} for a new $\gamma$}. In our examples, we only call this with two values of $\gamma$.   

 \item $O(N)$ cost to
form $\bfW$ for each time the algorithm is called.   
\end{enumerate}
}

\subsection{Adaptations for Shift Range} 

In our application, as in many others where our two-level recycling strategy may be employed, the shifts span many orders of magnitude.  We may need to adjust the two level approach to account for the wide range of shifts.  We do this by partitioning the shifts into groups by magnitude.
We  use a partition into {\bf two} groups, though we note that the approach we now discuss can be easily adapted to more than two clusters of indices.  The notation is as follows: 
\[  \mbox{Left group: shift indices } 1 \le \ell \le \sI \qquad \mbox{Right group: shift indicies } \sI < \ell \le M. \]
In our application, the $\lambda_i$ are log-evenly spaced and range over
at least five orders of magnitude, so the $\gamma_i$ span at least ten orders of magnitude.  Thus, the grouping is not even, but rather will be guided by the analysis in Section \ref{sec:analysis}.  For example, for 20 log-evenly space $\lambda_i$, the left grouping might contain 15 of the shifts and the right 5.

\subsubsection{Local Recycle spaces}
In our application, since the $\gamma_\ell$ vary significantly,
the local shift-specific subspaces will be seeded differently depending on if they are in the right or the left grouping of 
indices.  Algorithm (\ref{alg:updatelocalV2}) gives the modification to Algorithm (\ref{alg:updatelocal}) in this case.    

\begin{algorithm}[ht] \caption{ \label{alg:updatelocalV2} 
 Update left and right shift-specific recycle spaces $\bfUU_{L,\ell}$ or $\bfUU_{R,\ell}$ }
 Given $1 \le J_L \le \mathcal{I}$ (or $\mathcal{I} < J_R \le M)$ \\
 Call Algorithm 1 with $\gamma=\gamma_{J_L}$ (or $\gamma = \gamma_{J_R})$ and return $\MW$ \\
 Update $\bfUU_{L,\ell}$ (or $\bfUU_{R,\ell})$ as $[\MW,\bfg^{(k-1,\ell)},\bfg^{(k,\ell-1)}]$ \\
 \end{algorithm}
  
Thus, the updates for the current $k$ over all $\ell$ are given by 
\begin{equation}
\Xkl = \left\{ \begin{array}{ll}  \bfy_{m_{\ell}}^{\ell} + \bfUU_{L,\ell} \bfK_{L,\ell}^T [ \bfrkl - \left( \bfA + \gamma_\ell \bfE_k \right) \bfy_{m_{\ell}}^{\ell} ] + \Xkl_0 & \ell \le \sI \\
      \bfy_{m_{\ell}}^{\ell} + \bfUU_{R,\ell} \bfK_{R,\ell}^T [ \bfrkl - \left( \bfA + \gamma_\ell \bfE_k \right) \bfy_{m_{\ell}}^{\ell}] + \Xkl_0 & \ell > \sI \end{array} \right. .
\label{eq:upest} 
\end{equation}

\subsection{Algorithm}  The complete algorithm is outlined in Algorithm \ref{alg:shifted}. For reference, the index notation is summarized in Table \ref{tab:indicies}.  Note that, for updates to the local recycle spaces, we have assumed 
a single partitioning of the shifts into a right and a left region as discussed earlier, but if warranted, there could be further subdivision.  

\begin{table}[h]
    \centering
    \begin{tabular}{c|l} \hline
        $k$  & index of the outer iteration \\ \hline
        $\ell$ & indexes the shifts ($\gamma=\gamma_\ell$) \\ \hline
          $M$   & total number of shifts \\ \hline
     $\mathcal{I}$     & index between 1 and $M$ determining shift groupings \\ \hline
        $i_1$ & index corresponding to the first initial system solve (outer loop)\\ \hline
         $i_2$ & index corresponding to the second initial system solve (outer loop) \\ \hline
         $i_*$ & index of system to establish principle $\MU$ and $\MK$ (\ref{eq:establishU}). \\ \hline
         $\ell_c$ & index of correction system (\ref{eq:delta}) solved to update the principle recycle space \\ \hline
         $J_R$ & index on right for generating a local recycle space via Algorithm 4 \\ \hline
         $J_L$ & index on left for generating a local recycle space via Algorithm 4 \\ \hline
    \end{tabular}
    \caption{Summary of index notation in the main algorithm.}
    \label{tab:indicies}
\end{table}

\begin{algorithm}[t]  \caption{\label{alg:shifted} Shifted Recycling Algorithm}
\SetAlgoLined
Given $\gamma_1,\ldots,\gamma_M$, and index $\mathcal{I}$ with $1 \le \mathcal{I} \le M$ to divide the left and right regions. \\ 
Solve  $\bfx^{(0,i_1)}, \bfx^{(0,i_2)}$ via MINRES to tolerance. \\
Use MINRES output to approx. eigenspaces $\bfV_{i_1}, \bfV_{i_2}$ of 2 initial systems \\
\For{$k=0$ to maxits}{
  \If{k == 0}{ 
 $\widetilde{\bfUU} = [\bfx^{(0,i_1)},\bfx^{(0,i_2)},\bfV_{i_1},\bfV_{i_2}]$,}
\Else{ 
 Compute $\bfV^{(new)}$ as described in Section (\ref{ssec:updatemaster}). \\
$\widetilde{\bfUU} = [\bfV_1,\bfV^{(new)},\bfx^{(k-1,1)},\ldots,\bfx^{(k-1,M)},]$}
Use Algorithm \ref{alg:orthU} to update $\widetilde{\bfUU}$. \\
  Choose $\gamma_*$ to be one of the shifts; compute $\bfUU$ (implicitly) and $\bfK$ from (\ref{eq:establishU}). \\
  Determine $\bfx^{(0,\ell)} = \bfUU \bfq$ where $\bfq$ is from (\ref{eq:orthupdate}). \\
   Determine $\hat{\bfA}, \hat{\bfE}$ from (\ref{eq:hq}). \\
    \For{$\ell = 1:\mathcal{I}$}{
   Call Algorithm \ref{alg:updatelocalV2} with $J_L$ to find $\bfUU_{L,\ell}$.
Run MINRES for $m_{\ell}$ steps to get the updated estimate $\bfx^{(k,\ell)}$ via  (\ref{eq:upest})}  %close the for ell = .. I
\For{$\ell=\mathcal{I}+1:M$}{
Call Algorithm \ref{alg:updatelocalV2} with $J_R$ to find $\bfUU_{R,\ell}$. 
Run MINRES for $m_{\ell}$ steps to get the updated estimate $\bfx^{(k,\ell)}$ via (\ref{eq:upest})
} %close for ell
Compute L-curve and pick $\ell_*$ \\
From $\bfx^{k,\ell_*}$, generate $\bfD^{(k+1)}$ \\
$\bfE \leftarrow \bfL^T 
\left(\bfD^{(k+1)}\right)^2 \bfL$ \\
}
\end{algorithm}

\section{Overhead Cost Summary}
\label{sec:overhead}

While the dominant costs remain the total number of matrix-vector products to solve the systems for all shifts, we briefly
outline the overhead associated with 
Algorithm \ref{alg:shifted}.  
Specifically, we summarize the 
$O(N)$ or greater effort involved in establishing subspaces to recycle, in computing initial guesses, etc.  

\medskip
\noindent{\bf Initialize $\widetilde{\bfUU}$ (Steps 2-3,6)}.  With the calculation of the initial space $\widetilde{\bfUU}$ we need
 \begin{enumerate}
 \item $m_{i_1} + m_{i_2}$ matvecs with $\MA$ and $\ME$ to solve the two systems to tolerance; 
  \item if $K$ Ritz vectors are desired for each of the 2 systems, we also need orthogonalization of $2K$ Lanczos vectors, $2K$ matvecs with $\MA, \ME$ to form the 2 projected systems, and $O(N K^2)$ flops to form $\MV_{i_1},\MV_{i_2}$.
  \end{enumerate}

\medskip
\noindent{\bf Update Principle Recycle Space (Steps 9-10)}:  Updating the principle recycle space requires computing $\MV^{(new)}$ with a maximum cost of 100 matvecs with $\MA, \ME$.  
  
\medskip
\noindent{\bf Find well-conditioned $\widetilde{\bfUU}$ (Step 12)} This step requires an SVD for a matrix with $N$ rows and $O(100)$ columns.

\medskip
\noindent{\bf Finding $x^{(0,\ell)}$ (Step 14):} When $\widetilde{\bfUU}$ is known, 
overhead in
finding $x^{(0,\ell)}$ includes:
\begin{enumerate}
    \item $n_c$ matvecs with $\MA, \ME$ in (\ref{eq:establishU});
    \item orthogonalization of $n_c$ vectors to get $\MK$ in (\ref{eq:establishU});
\end{enumerate}
Note that $\MU$ does not need to be computed explicitly (avoiding $O(N n_c^2)$ work). Multiplication by $\MU$ (or by $\ME_k \MU$) can be done by first computing, say, $\M{R}^{-1}\Vq$ and then computing $\widetilde{\MU} (\M{R}^{-1}\Vq)$. Analogously, $\ME_k \MU\Vq = \ME_k \widetilde{\MU}(\M{R}^{-1}\Vq)$.
The matrices $\MK^T (\ME_k \MU)$, $(\ME_k \MU)^T(\ME_k \MU)$, and $(\MU_k^T\ME_k)\V{b}$ can all be computed from $\ME_k \widetilde{\MU}$ and 
subsequent multiplication of $n_c \times n_c$ matrices with $\M{R}^{-1}$ without additional $O(N)$ cost.

\medskip
\noindent{\bf Updating local recycle spaces (Step 15,17,20)}.  Here, we incur
\begin{enumerate}
 \item the cost to compute each of the projections $\hat{\bfA}, \hat{\bfE}_k$, which {\it only needs to be done once per outer iteration $k$, independent of how many calls to Algorithm \ref{alg:interior} for a new $\gamma$};    
 \item an $O(N)$ cost to
form $\bfW$ for each time the algorithm is called (twice, in our case);
 \item and an $O(MN)$ orthogonalization cost to find the $\bfK_{L,\ell}$ $(\bfK_{R,\ell})$.
\end{enumerate}

\medskip
Note that the extra matvecs in initializing and updating the principle recycle space are accounted for in graphs in the numerical results, while the other overhead counts are not.  Nevertheless, as the numerical results reveal, with the savings in matvecs being on the order of 10,000, the remaining overhead is marginal in comparison.  
        
\section{Alternative Strategy for Generating Starting Guesses}
\label{sec:oblique}

As our method is more broadly applicable to sequences of shifted systems,  we briefly discuss a potentially cheaper strategy for starting guesses.  
Given the initial residual 
\eqref{eq:res_initGuessPrincRecy}, 
we can compute $\Vq$,  
for each 
$\d_\ell = \c_\ell - \c_{*}$, using the oblique projection 
$\MK^T \Vr_0^{(k,\ell)} = \V0$:
\eqs
  \MK^T \Vb & = & (\MI + \d_\ell \MK^T\ME_k\MU)\Vq ,
  \label{eq:obliqueupdate}
\eqe
which still requires the one time computation
of $\MK^T \ME_{k} \MU$ but saves the 
computation of $(\ME_k\MU)^T(\ME_k\MU)$. 
In this case, 
we need to solve (\ref{eq:obliqueupdate}) 
for each shift, and compute $\bfx_{0}^{(k,\ell)}$ accordingly.  
Using the oblique projection updates {\it for all $\ell$} saves $O(N n_c^2)$ flops over the orthogonal projection approach.

\subsection{Analysis: Oblique vs. Orthogonal Projection}  
We briefly analyze the difference between the orthogonal and oblique projections and the resulting 
reductions
in the initial residual norm.
% relative to $\| \Vb \|_2$.
For this analysis, let the columns of $\MK$, which defines the oblique projection, be  orthonormal, let $\MY \MS$ be the economy QR decomposition of $\MK + \d_\ell \ME_k\MU$ 
($\MY \in \Rmn{N}{n_c}$), 
and assume that 
$\MK^T\MY$ is invertible. Let $\Vq_1$
(optimal) and $\Vq_2$ be given by the projections 
\eqs \label{eq:orthoP}
  \Vb - \MY \Vq_1 & \perp & \MY  \quad \eqv \quad 
  \Vq_1 = \MY^T \Vb , 
\\
\label{eq:oblP}
  \Vb - \MY \Vq_2 & \perp & \MK \quad \eqv \quad 
  \Vq_2 = (\MK^T \MY)^{-1}\MK^T \Vb , 
\eqe
with best approximation $\MY \Vq_1 = \MY \MY^T \Vb$ and residual $\Vr_1 = (\MI - \MY \MY^T) \Vb$, and 
`oblique' approximation
$\MY \Vq_2 = \MY(\MK^T \MY)^{-1}\MK^T \Vb$
and residual 
$\Vr_2 = (\MI - \MY(\MK^T \MY)^{-1}\MK^T)\Vb$.
The corresponding, respective, initial guesses are $\MU \MS^{-1}\Vq_1$ and 
$\MU \MS^{-1}\Vq_2$.
Using the singular value decomposition (SVD) $\MY^T\MK = \MgF \MgO \MgY^T$, 
we can give simple expressions for 
$\Vq_1 - \Vq_2$, 
the difference between the oblique and orthogonal projection, $\MY(\Vq_1 - \Vq_2)$, and the difference between the residuals, $\Vr_1 - \Vr_2$. 
Given that 
$\MY^T\MK$ is invertible, the singular values $\o_j > 0$.
From this SVD, we have 
$\MY \MY^T \MK\MgY = \MY \MgF \MgO$
and
\eqs \label{eq:varCS}
  (\MI  - \MY\MY^T)\MK\MgY = 
  \MY_{c} \MgSi 
  & \quad\eqv\quad &
  \MK\MgY = \MY \MgF \MgO + \MY_c \MgSi ,
\eqe
where
$\MgO^2 + \MgSi^2 = \MI$,
the nonzero columns of $\MY_c$ are orthonormal\footnote{The columns of 
$(\MI  - \MY\MY^T)\MK\MgY$ are
orthogonal, as can be verified 
directly by taking their inner products. Alternatively, note that 
\eqref{eq:varCS} is part of the 
generalized CS-decomposition
in \cite{PaigSaun_1981}.}, 
and by construction $\MY^T \MY_c = \M0$.
If $\o_j = 1$, then 
$\s_j = 0$, and the corresponding column in $\MY_c$ can be chosen arbitrarily; for convenience,
here, we choose $(\MY_c)_j = \V0$.
Using the SVD of $\MK^T\MY$ in 
\eqref{eq:oblP} gives 
$\Vq_2 = \MgF \MgO^{-1}\MgY^T\MK^T\Vb$,
and substituting for 
$(\MK\MgY)^T$ using \eqref{eq:varCS} gives
\eqs
  \Vq_2 =
  \MgF \MgO^{-1}
  [\MgO\MgF^T\MY^T + \MgSi \MY_c^T]\Vb 
  & = & 
  \MY^T \Vb + 
  \MgF\MgO^{-1}\MgSi\MY_c^T\Vb.
\eqe
As $\Vq_1 = \MY^T\Vb$, 
$\Vq_2 - \Vq_1 = 
\MgF\MgO^{-1}\MgSi\MY_c^T\Vb$,
and the difference between the
projections (and the residuals) is given by 
\eqs
  \MY(\Vq_2 - \Vq_1) 
  = \Vr_1 - \Vr_2 \;\;
  & = &
  \sum_{\{j \;|\; \o_j < 1\}} 
  (\MY\Vgf_j) \frac{\s_j}{\o_j} (\MY_c)_j^T \Vb.
\eqe
Since, by construction $\Vr_1 \perp \MY$, we have $\|\Vr_2\|_2^2 = \|\Vr_1\|_2^2 + 
\| \MgO^{-1}\MgSi\MY_c^T\Vb \|_2^2$.
Notice that if $\Vb \in \Ra{\MY}$, then 
$\Vq_2 = \Vq_1$ and $\Vr_2 = \Vr_1 = \V0$.
We see that $\| \Vr_2 \|_2 \gg \|\Vr_1\|_2$ can only happen if 
at least one of the $\o_j$ is small
and the corresponding
$(\MY_c)_j^T \Vb$ is (relatively) large.
Moreover, $\| \MY_c^T \Vb \|_2
\leq \|(\MI-\MY\MY^T)\Vb\|_2 = \|\Vr_1\|_2$; so, 
very small $\|\Vr_1\|_2$ implies
very small $\| \MY_c^T \Vb \|_2$
and hence generally $\|\Vr_2\|_2$ 
is small as well, unless some of the $\o_j$ are very small.
This gives additional motivation to seed the recycle space with approximate solutions, especially when an oblique projection is competitive.

\subsection{Implementation Considerations} \label{sssec:obliqueinit}

Recall $\bfUU$ was obtained from a single choice in shift in (\ref{eq:establishU}).   If the range in the shifts is large enough, we do not expect the oblique projection update for computing the true solution to be suitable across all shifts $\gamma_\ell$. In this case, we can modify our approach. 
Let $I_L$ denote an index $1 \le I_L \le \sI$ (that is, an index within the left group) and $I_R$ denote an index in the right group, such that
$\sI < I_R \le M$.   We first compute 
\[  (\bfA + \gamma_{I_L} \bfE_k) \widetilde{\bfUU} = \bfK_L \bfR_L, \qquad 
    (\bfA + \gamma_{I_R} \bfE_k) \widetilde{\bfUU} = \bfK_R \bfR_R, \]
so that $\bfUU_L = \widetilde{\bfUU} \bfR_L^{-1}$, $\bfUU_R = \widetilde{\bfUU} \bfR_{R}^{-1}$.   Although the ranges of $\bfUU_L$ and $\bfUU_R$ are the same, the ranges of $\bfK_L, \bfK_R$ are not.  

This induces two different Petrov-Galerkin systems and solution updates.  However, otherwise we follow the approach, replacing $\bfUU, \bfK$ with 
$\bfUU_L$, $\bfK_L$ for the first (left) group and with $\bfUU_R,\bfK_R$
for the second (right) group. 
We now have
  \begin{equation} \label{eq:updates}
   \Xkl_0 := \left\{ \begin{array}{ll} \bfUU_L \left( \bfI + (\gamma_\ell - \gamma_{I_L}) \bfK_L^T \bfE_k \bfUU_L \right)^{-1} \bfK_L^T \bfb & \ell \le \sI \\
   \bfUU_R \left( \bfI + (\gamma_\ell - \gamma_{I_R}) \bfK_R^T \bfE_k \bfUU_R \right)^{-1}\bfK_R^T \bfb  & \ell > \sI \end{array} \right.  .    
\end{equation}
However, with this modification, using  two systems, {\it we have extra overhead} in computing both $\bfK_L$ and $\bfK_R$ (a second reorthogonalization of an $N \times n_c$ matrix is required) and in computing both $\bfK_L^T \bfE_k \bfUU_L$ and $\bfK_R^T \bfE_k \bfUU_R$ one time for all the indices.  This means that
the {\it cost savings over the orthogonal projection approach has effectively disappeared, } suggesting  that in applications where the shifts vary widely, such as the applications in our paper, the orthogonal projection approach is to be preferred.

\section{Numerical Results} \label{sec:num}

The numerical results are all computed using Matlab 2022a.  First, we need to discuss the selection of the stopping criterion 
in the context of the imaging problem and computing the regularization parameter at each stage.

\subsection{Residual Tolerance} \label{ssec:resi}

MINRES itself is known to act as a regularization method \cite{Hanke,KilmerStewart}.  Thus, if we stop iterating on (\ref{eq:shift}) prior to reaching the convergence tolerance on the residual, we expect the solution to be smooth and not contaminated by noise.    This is a problem for the small shifts in our context, because we need to 
see the impact of $\lambda$ as the regularization parameter. Therefore, 
we have to solve each new system for $\lambda$ accurately enough. If we stop too soon, we only see the impact of the iteration index as the  regularization parameter,   
and points corresponding to small $\lambda$ will be effectively indistinguishable on the L-curve. This makes it difficult to determine a good value of the regularization parameter, because the vertical arm of the L-curve is almost non-existent.   
Thus, we propose to solve each system so that the relative residual norm of the regularized normal equations system is $10^{-6}$.

\subsection{Example 1:  Image Deblurring} \label{ssec:one}

\begin{figure}
\centering
\subfigure[Original.]{\label{fig:a1} \includegraphics[height=1.25in,width=1.75in]{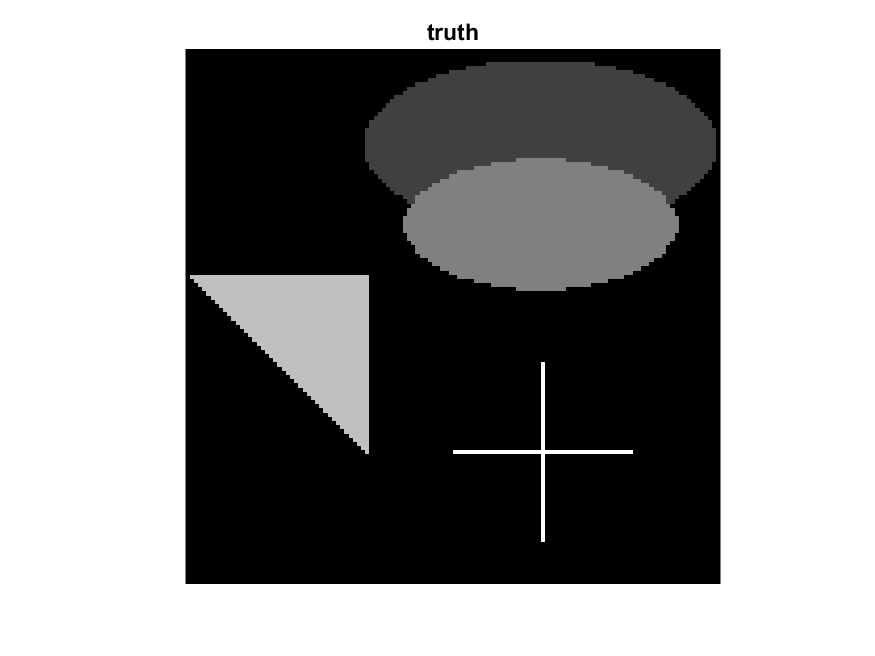}}
\subfigure[Blurred + Noisy.]{\label{fig:a2} \includegraphics[height=1.25in,width=1.75in]{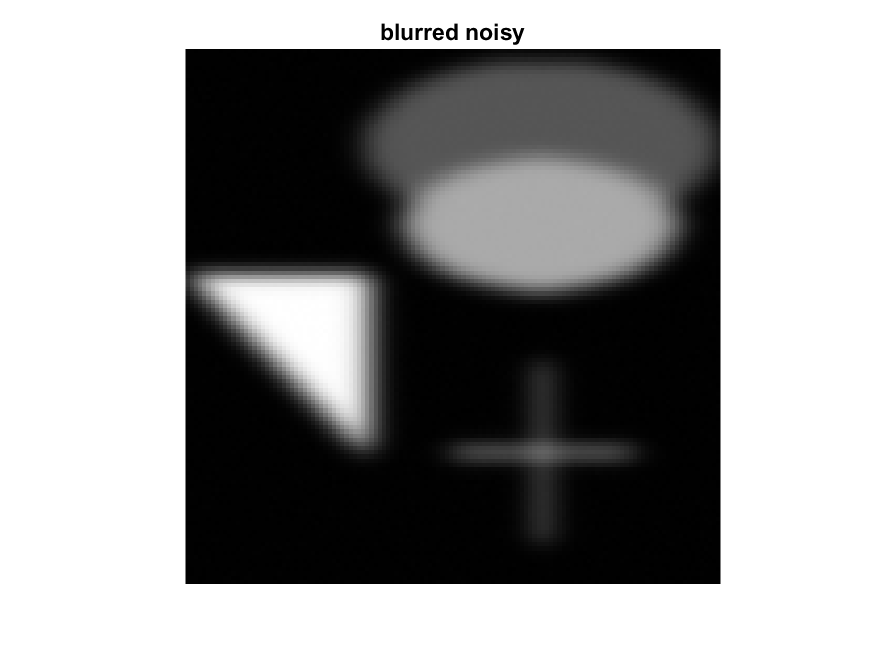}}
\subfigure[First Recon.]{\label{fig:a3} \includegraphics[height=1.25in,width=1.75in]{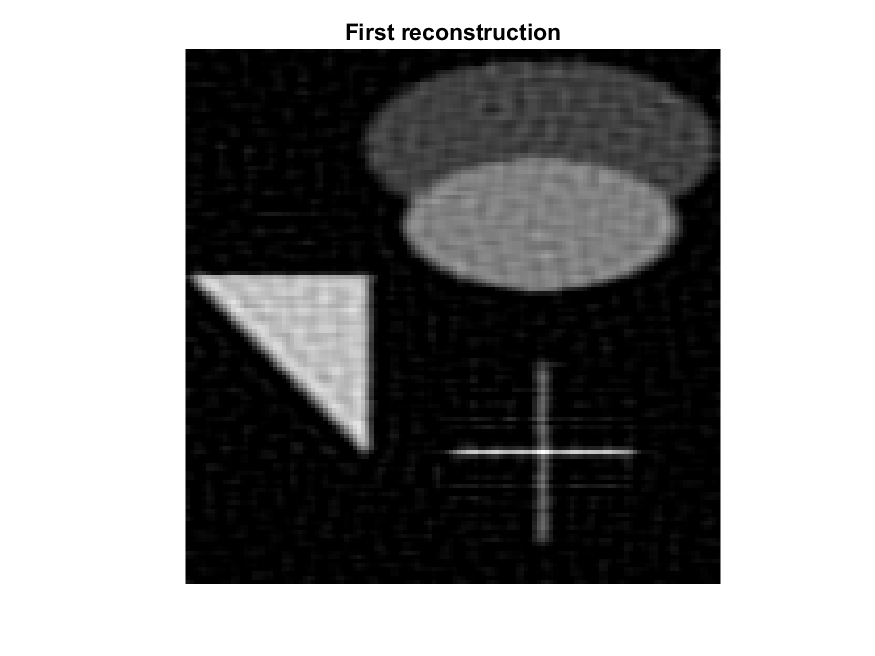}}
\subfigure[Final Recon.]{\label{fig:a4} \includegraphics[height=1.25in,width=1.75in]{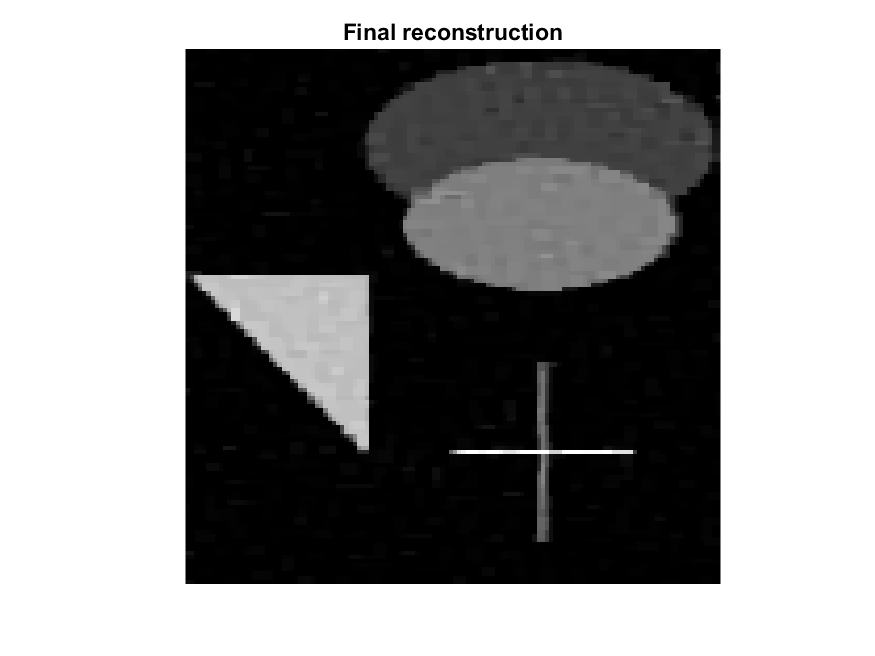}}
\caption{ \label{figex1:images} Comparison of the images in Example 1.}
\end{figure}

We consider a symmetric, Gaussian blur and zero boundary conditions.  In this example,  
\[  \bfC := \sum_{i=1}^2 C_i^{(1)} \otimes C_{i}^{(2)}, \]
with each $C_{(i)}^{(j)}$ a Toeplitz matrix of differing bandwidth and 
$\sigma$. Here, the bandwidths for $C_1^{(1)}, C_1^{(2)}, C_2^{(1)}, C_2^{(2)}$ were 8,7,12 and 4, respectively; the $\sigma$ were 3,$\frac{5}{2}$,4 and 2, respectively. 
 The ground truth $128 \times 128$ image in Figure \ref{figex1:images}(a) is taken from the example {\tt blur.m}  in \cite{Regtools}.  
The blurred image was obtained as the product $\bfC \mbox{ vec}(\bfX_{true})$, and the data was determined by adding 
0.5 percent Gaussian noise to the result.  The blurred, noisy data image is displayed in Figure \ref{figex1:images}(b).

Here, $\lambda$ was determined from 20, log-evenly spaced values between $10^{-4}$ and $10^{2.5}$, so $M=20$.  
In this example, we determine the left set of systems to correspond to the smallest 15 values ($\mathcal{I}=15$) and the right as the remaining five systems.
We use $i_1=1$ (i.e., the smallest possible shift) and $i_2 = 10$.  We used the first 100 Krylov vectors (reorthogonalized) of each run to determine the Ritz vectors in $\MV_{i_1}$ and $\MV_{i_2}$, respectively.  We kept all 100 eigenvector approximations and placed them in $\MV_{i_1}$ whereas we kept only 50 eigenvector approximations in $\MV_{i_2}$, so that the dimension of the initial space $\widetilde{\MU}$ is 153 (see \ref{eq:widetildeU}).   
We take $i_* = 10$,    
$J_L=15$ and $J_R=19$.   The number of local Ritz vectors we take in Algorithm \ref{alg:updatelocal} is 12, and thus 
the maximum dimension of any of the local recycle spaces $\MU_\ell$ is 14.  The index $\ell_c$ of the correction system is taken to be 18.

For terminating the outer iterations, we use the stopping criteria from \cite{MaclachlanBoltenKilmer2023}: once the selected $\lambda$ is the same for three consecutive iterations, we stop.   In this example, the outer iteration terminates after 11 iterations.   
The final reconstruction is 
given in Figure \ref{figex1:images}(d).

To minimize the overhead in computing $\bfK_{L,\ell}$ (and similarly for $\bfK_{R,\ell}$), since the $\widetilde{\bfUU}_{L,\ell}$ ($\bfUU_{R,\ell}$) share columns as a function of $\ell$, we need only compute the corresponding matvecs once.  Furthermore, the orthogonalization of the first 12 columns in computing the $\MK_{\ell}$ does not have to be repeated, either, as a function of $\ell$.  

%\comment{
\begin{figure}
    \centering
    \includegraphics[width=3in,height=2in]{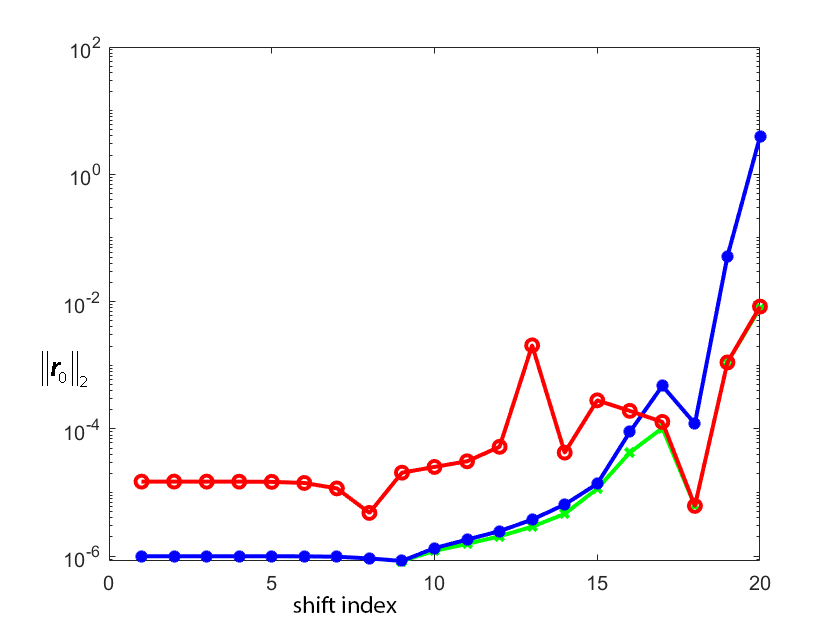}
    \caption{Although the orthogonal projection method in Section \ref{ssec:initorth} is used in this experiment to generate initial guesses for each outer iteration $k$, we demonstrate here the alternative oblique projection strategy for outer iteration $k=6$. The green curve corresponds to the optimal, but too expensive, case of recomputing $\bfK$ for every shift and explicitly doing orthogonal projection (in which case $\bfq = \bfK^T \bfb$).  
    The blue curve shows residual for initial guesses via oblique projection fixing  $\lambda_*$ as a small shift (see (\ref{eq:obliqueupdate})); the red curve shows the result if a larger shift is used for the oblique projection.  Note the blue curve is almost optimal 
    through index 15, while the red curve overlays the optimal curve for indices 17 and higher.
    In practice, we would use the small shift for $i \le \mathcal{I}=15$ and the large shift for $i > \mathcal{I}$ 
     as in (\ref{eq:updates}).  }
    \label{fig:my_label}
\end{figure}

%\subsubsection{Initial Solution Estimates}

\begin{figure}
\centering
\subfigure[Matvecs, No Recycle.]{\label{fig:b1} \includegraphics[scale=.4]{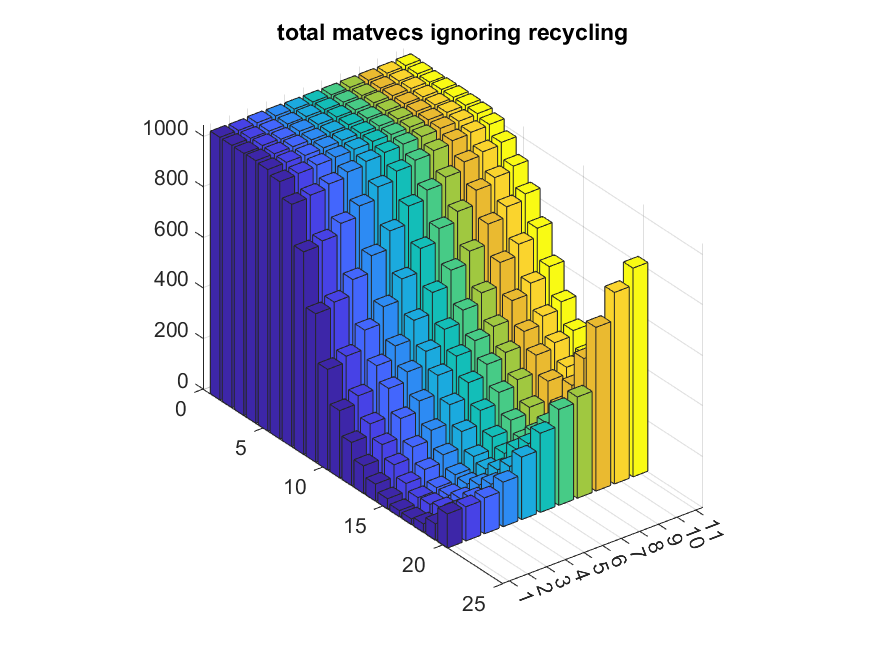}}
\subfigure[Matvecs with Recycling.]{\label{fig:b2} \includegraphics[scale=.4]{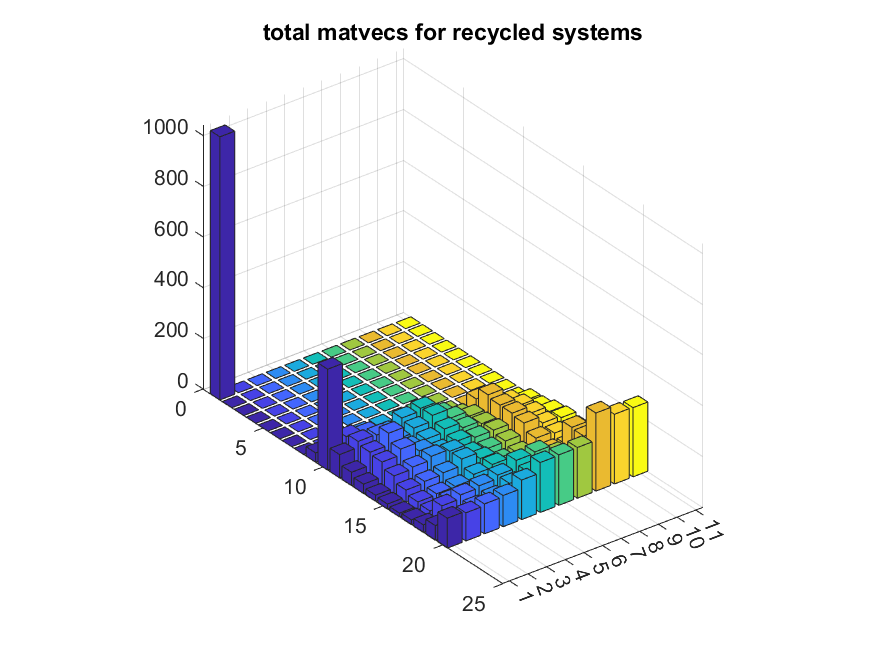}}
\subfigure[Savings in Matvecs.]{\label{fig:b3} \includegraphics[scale=.4]{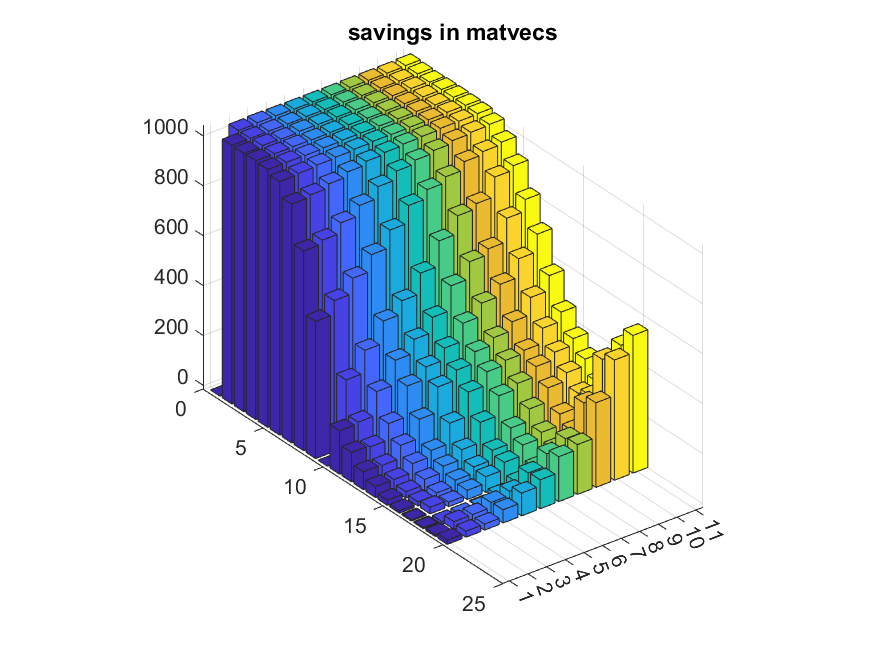}}
\caption{ \label{ex1:matvecs} Number of matvecs that would be needed to solve each of the current systems if recycling were not used; number of matvecs when recycling according to our method, and total savings in matvecs for Example 1.}
\end{figure}

As is clear from Figure \ref{ex1:matvecs}, the overall savings in matrix-vector products achieved by our new approach vs. naively solving each linear system to tolerance at each outer iteration is significant.  While there is some overhead in terms of additional matvecs in computing initial guesses (which are accounted for in the plots), some overhead in updating the principle recycle space (also accounted for in the plots),  and overhead updating the local recycle spaces, these costs are clearly more than amortized by the overall savings in matvecs.  The spikes in  matvecs for $k=1$ are due to the one-time cost of solving 2 shifted systems to tolerance.  The cost 
% with computing 
of solving (\ref{eq:delta}) to improve the principle space for the larger shift systems is included in the plots, but clearly does not constitute a large overhead. 

\subsection{Example 2: Computed Tomography}

\begin{figure}
\centering
\subfigure[Original.]{\label{fig:c1} \includegraphics[height=1.25in,width=1.75in]{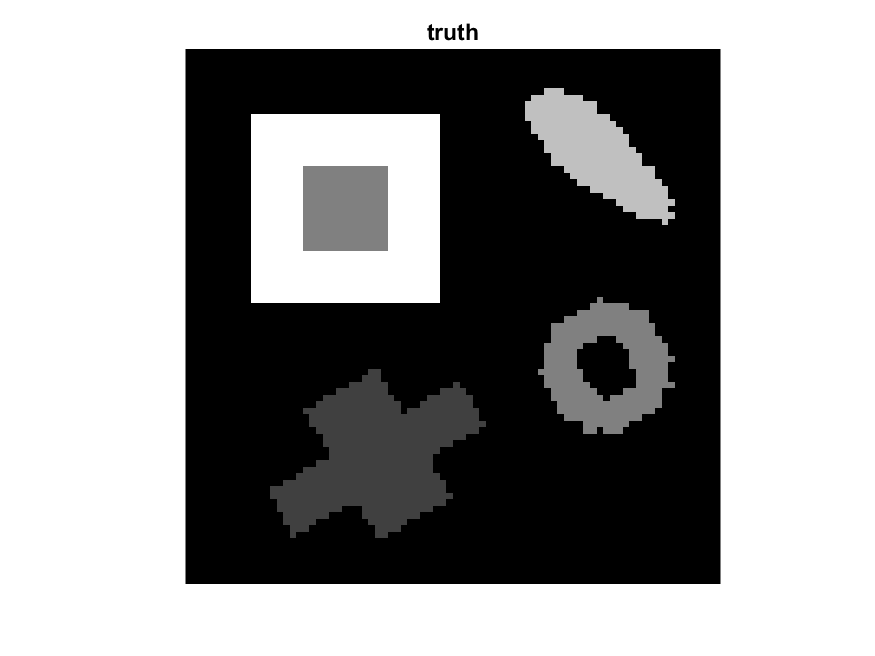}}
\subfigure[First Recon.]{\label{fig:c2} \includegraphics[height=1.25in,width=1.75in]{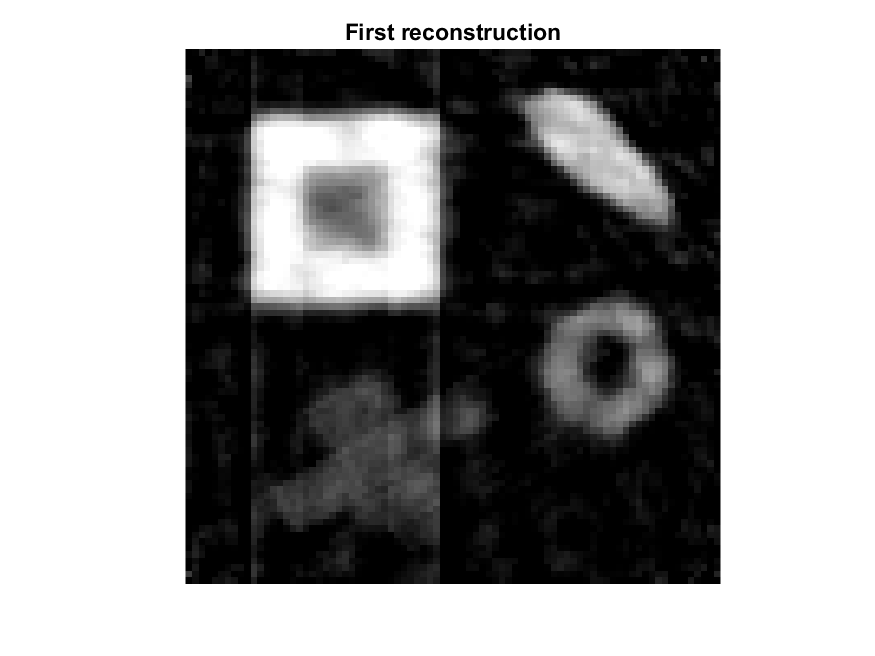}}
\subfigure[Final Recon.]{\label{fig:c3} \includegraphics[height=1.25in,width=1.75in]{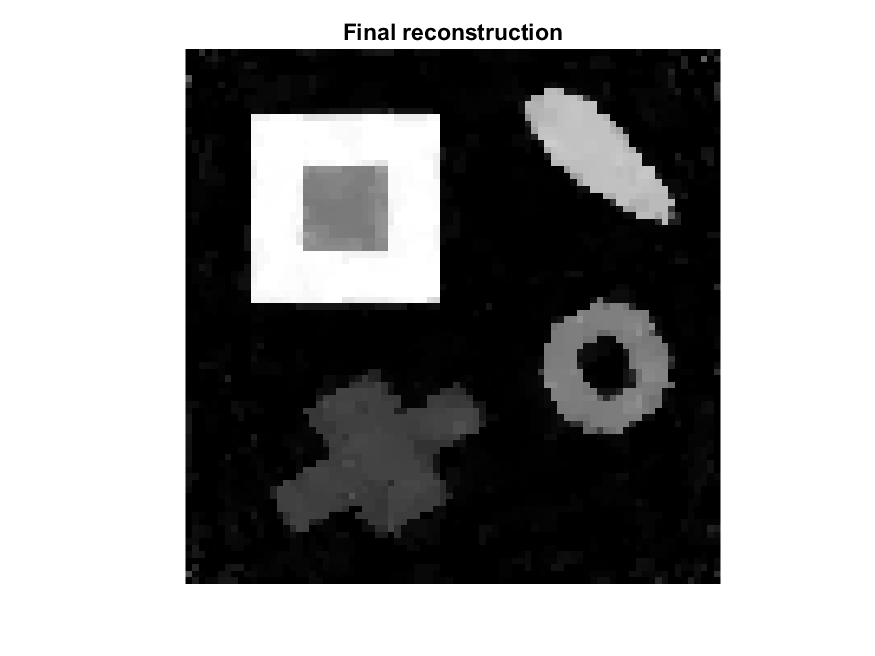}}
\caption{ \label{figex2:images} Comparison of the images in Example 2.}
\end{figure}

In this example, the true image is the $n \times n$ image in Figure \ref{fig:c1}.  Here, $n=82$.  The operator $\M{C}$ is the discrete Radon transform matrix for parallel beam CT corresponding to taking measurements angles from 1 to 135 degrees in increments of 1 degree.\footnote{The matrix can be constructed in Matlab column by column by using as input the standard unit vectors $\V{e}_i$, $i=1,\ldots,n^2$.}  We normalize so that $\M{C}$ has unit Frobenius norm. We form the true ``sinogram data" as the product $\M{C} \V{x}_{true}$ and add 1 percent Gaussian white noise to simulate the measured data, as in (\ref{eq:forward}).  

For this problem, we use 20 equally spaced values ($M = 20$) of $\lambda$ between $10^{-4}$ and $10^{1}$.  
We again use $i_1 = 1$ and $i_2 = 10$ for 
the indices of the systems whose solutions we generate outside the loop, and whose Krylov vectors we use 
%to make the subspace we use 
to initialize the principle recycle space. 
We take $\gamma_*$ as the 12th smallest shift. 
We choose $\mathcal{I} = 16$, so the left subgroup of indices contains 16 systems, the remaining 4 systems belong to the right subgroup.   The index that defines the correction system (5.2) is $\ell_c = 19$, 
and $J_L$ ($J_R$) that specifies which 
system's eigenvectors to approximate is 15 (19).

\begin{figure}
\centering
\subfigure[Matvecs, No Recycle.]{\label{fig:d1} \includegraphics[scale=.4]{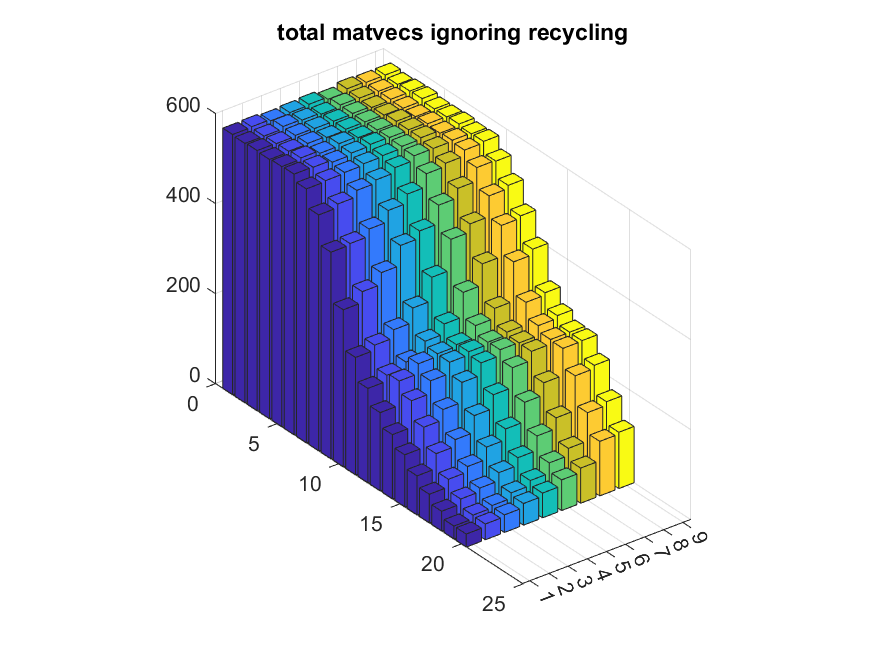}}
\subfigure[Matvecs with Recycling.]{\label{fig:d2} \includegraphics[scale=.4]{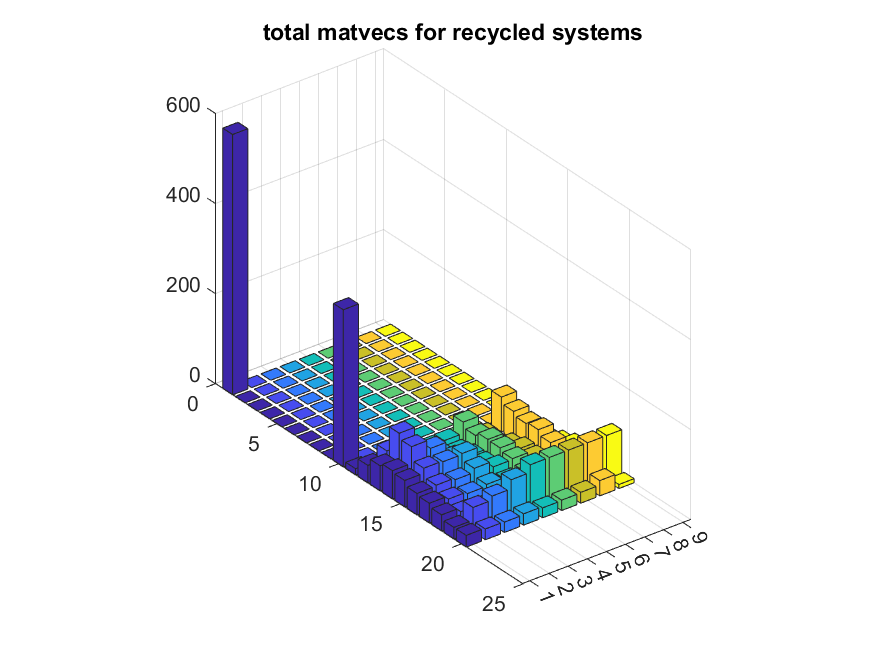}}
\subfigure[Savings in Matvecs.]{\label{fig:d3} \includegraphics[scale=.4]{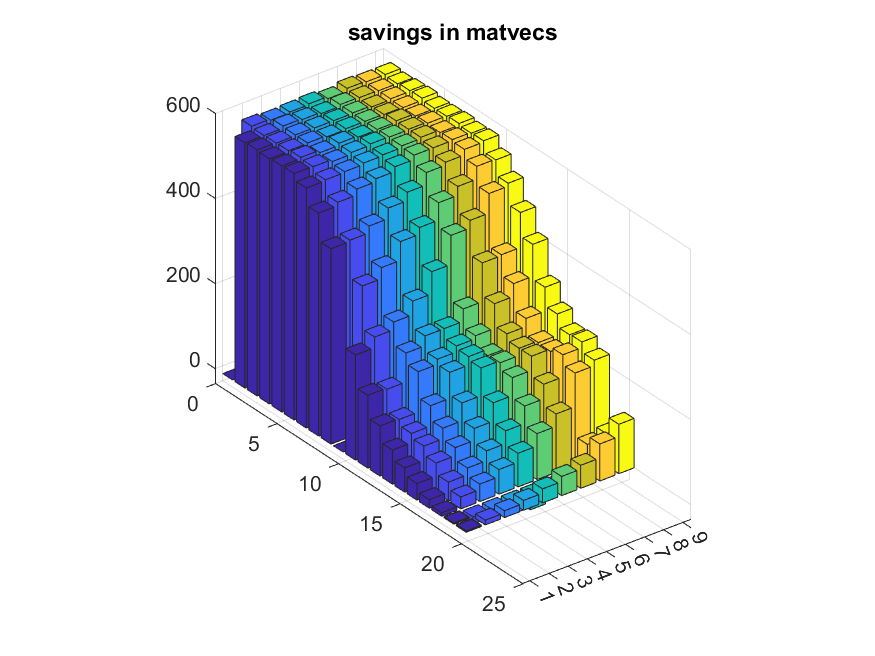}}
\caption{\label{ex2:matvecs}Number of matvecs to solve each of the current systems without recycling; number of matvecs when recycling according to our method, and total savings in matvecs for Example 2.}
\end{figure}

Using the stopping criteria in \cite{MaclachlanBoltenKilmer2023} for the outer iteration leads to the algorithm taking a total of 9 outer iterations.  Figure \ref{ex2:matvecs} illustrates the substantial savings in matvecs across all shifted systems and all outer iterations.  As the outer iteration increases and edges are picked up, the matrix $\bfE_k$ becomes more ill-conditioned, resulting in an increased number of matvecs required to solve the systems with the larger shifts.  However, the recycling strategy we employ here clearly mitigates this effect, keeping the number of required matvecs quite low.

\section{Conclusions and Future Work}  \label{sec:conclusions}

In this paper, we introduce a new inner-outer Krylov subspace recycling approach for large-scale image reconstruction and restoration techniques that employs a nonlinear iteration to compute 
%the right 
a suitable
regularization matrix. By maintaining a larger dimensional principle recycle space that is rich in information shared across systems with varying shifts, we can build good initial solutions after each nonlinear update of the regularization matrix. To solve the individual correction equations, we employ recycling of local subspaces, which are informed by the principle recycle space as well as local information and prior solutions.  

We provide an analysis based on the  generalized eigenvectors to motivate reusing solutions 
to seed recycle spaces.  The cost analysis and numerical results show the huge potential savings over naively solving the systems individually with MINRES.

Though motivated by the application of edge preserving regularization, our techniques are suitable in other applications involving non-identity shifted systems with symmetric and positive (semi)definite matrices.  To this end, we also give an alternative means for computing initial guesses based on oblique projection.  This would be cheaper for problems where problem characteristics do not change as dramatically across shifts as they do in the present application.

Finally, we note that all the iterations were unpreconditioned.  Recent work in \cite{MaclachlanBoltenKilmer2023} shows, for the particular application at hand, the potential benefit of additionally employing multigrid preconditioning for the systems dominated by the diffusion-like term in the problem.  In future work, we will consider how to combine the best features of the recycling approach presented here with the benefits of preconditioning for some, or all, of the shifted systems.  

\section{Acknowledgements}
The authors would like to thank Dr. Meghan J. O'Connell for technical discussions that led to the start of this paper.”

\bibliographystyle{abbrv}
\bibliography{Bibliography}
\end{document}

%% file: rom_macros.tex
\usepackage[normalem]{ulem}
\usepackage{amsmath}
\usepackage{amsfonts}
\usepackage{amssymb}
\usepackage{graphicx}
\usepackage{epstopdf,euscript}
\usepackage{url}
\usepackage{multirow}
\usepackage{array}
\usepackage{color}
\usepackage[ruled, linesnumbered]{algorithm2e}
\usepackage{mathtools}
\usepackage{hyperref}

\newcommand{\Xkl}{\mathbf{x}^{(k,\ell)}}

\newcommand{\bfrkl}{\mathbf{r}^{(k,\ell)}}

\newcommand{\bfd}{\mathbf{d}}
\def\bq{\begin{quotation}}
\def\eq{\end{quotation}}

\def\c{\gamma}

\def\d{\delta}

\def\m{\mu}

\def\s{\sigma}
\def\Si{\Sigma} % for SIAM Jrns \S is section symbol

\def\f{\varphi}

\def\F{\Phi}

\def\Y{\Psi}
\def\o{\omega}
\def\O{\Omega}

\newcommand{\MgF}{\boldsymbol{\F}}

\newcommand{\MgSi}{\boldsymbol{\Si}}
\newcommand{\MgO}{\boldsymbol{\O}}
\newcommand{\MgY}{\boldsymbol{\Y}}

\newcommand{\Vgf}{\boldsymbol{\f}}

\newcommand{\Rmn}[2]{\mathbb{R}^{#1 \times #2}}

\def\eqv{\Leftrightarrow}

\def\2nm#1{\|#1\|_2}

\def\Ra#1{\mathrm{Range}(#1)}

\newcommand{\ars}[1]{\left[ \begin{array}{#1}}
\newcommand{\are}{\end{array} \right] }
\newcommand{\oars}[1]{\begin{array}{#1}}
\newcommand{\oare}{\end{array}}
\newcommand{\rars}[1]{\left( \begin{array}{#1}}
\newcommand{\rare}{\end{array} \right) }

\newcommand{\eqs}{\begin{eqnarray}}
\newcommand{\eqe}{\end{eqnarray}}
\newcommand{\eqsn}{\begin{eqnarray*}}
\newcommand{\eqen}{\end{eqnarray*}}

\newcommand{\bmp}[2]{\begin{minipage}#1{#2}}
\newcommand{\emp}{\end{minipage}}

\newcommand{\ens}{\begin{enumerate}}
\newcommand{\ene}{\end{enumerate}}

\newcommand{\its}{\begin{itemize}}
\newcommand{\ite}{\end{itemize}}

\newcommand{\des}{\begin{description}}
\newcommand{\dee}{\end{description}}

\def\defs{\begin{definition}}
\def\defe{\end{definition}}
\def\teos{\begin{theorem}}
\def\teoe{\end{theorem}}
\def\prfs{\begin{proof}}
\def\prfe{\end{proof}}
\def\exas{\begin{exampl}}
\def\exae{\end{exampl}}
\def\excs{\begin{exercise}}
\def\exce{\end{exercise}}
\def\cors{\begin{corollary}}
\def\core{\end{corollary}}

\def\wtl{\widetilde}

\newcommand{\comment}[1]{} % makes its argument disappear

\newcommand{\bfDel}{\mbox{\boldmath$\Delta$}}

\newcommand{\bfx}{{\bf x}}
\newcommand{\bfK}{{\bf K}}

\newcommand{\bfy}{{\bf y}}

\newcommand{\bfv}{{\bf v}}
\newcommand{\bfb}{{\bf b}}

\newcommand{\bfr}{{\bf r}}

\newcommand{\bfz}{{\bf z}}
\newcommand{\bfq}{{\bf q}}
\newcommand{\bff}{{\bf f}}

\newcommand{\bfQ}{{\bf Q}}
\newcommand{\bfA}{{\bf A}}

\newcommand{\bfC}{{\bf C}}
\newcommand{\bfD}{{\bf D}}
\newcommand{\bfE}{{\bf E}}

\newcommand{\bfH}{{\bf H}}
\newcommand{\bfI}{{\bf I}}

\newcommand{\bfR}{{\bf R}}

\newcommand{\Co}{{\bf K}}

\newcommand{\bfX}{{\bf X}}

\newcommand{\bfUU}{{\bf U}}

\newcommand{\bfL}{{\bf L}}
\newcommand{\bfV}{{\bf V}}
\newcommand{\bfW}{{\bf W}}

\newcommand{\bfw}{{\bf w}}
\newcommand{\bfg}{{\bf g}}
\newcommand{\bfe}{{\bf e}}

\newcommand{\bea}{\left[ \begin{array} }
\newcommand{\eea}{ \end{array} \right] }

\newcommand{\mrg}{\mbox{Range}}